\documentclass[11pt]{amsart}

\usepackage[utf8]{inputenc}
\usepackage[a4paper,twoside]{geometry}

\usepackage[usenames,dvipsnames]{xcolor}
\usepackage{enumitem}
\usepackage{amsmath, amssymb, amsthm}

\usepackage[english]{babel}
\usepackage{fancyhdr}
\usepackage[all]{xy}

\usepackage{tikz}
\usetikzlibrary{calc,decorations.markings,patterns}

\usepackage[framemethod=tikz]{mdframed}

\usepackage{float}
\usepackage{multirow}
\usepackage[colorlinks=true,linkcolor=blue,citecolor=red,backref=page,pdftex]{hyperref}

\newcommand{\lib}{[}
\newcommand{\rib}{]}

\newcommand{\retrait}{\hspace{1.7cm}}

\newcommand{\unp}{\mathbf{\mathrm{1 \kern-0.25em I}}}
\newcommand{\un}{\mathbf{1}}

\newcommand{\R}{\mathbb R}
\newcommand{\N}{\mathbb N}
\newcommand{\C}{\mathbb C}

\newcommand{\Z}{\mathbb Z}

\newcommand{\di}{{\rm d}}

\newcommand{\cal}{\mathcal}

\newcommand{\supp}{\mathrm{supp}\,}

\newtheorem{theorem}{Theorem}[section]
\newtheorem*{theorem*}{Theorem}
\newtheorem*{lemma*}{Lemma}
\newtheorem*{corollary*}{Corollary}
\newtheorem{lemma}[theorem]{Lemma}
\newtheorem{proposition}[theorem]{Proposition}
\newtheorem*{proposition*}{Proposition}
\newtheorem{corollary}[theorem]{Corollary}

\theoremstyle{definition}
\newtheorem{definition}[theorem]{Definition}
\newtheorem*{definition*}{Definition}
\newtheorem{example}[theorem]{Example}

\newtheorem*{notations*}{Notations}
\theoremstyle{remark}
\newtheorem{remark}[theorem]{Remark}

\numberwithin{equation}{section}
\begin{document}

\title{The speed of convergence in the renewal theorem}

\author{Jean-Baptiste Boyer}
\email{jean-baptiste.boyer@math.u-bordeaux.fr}
\date{\today}

\begin{abstract}
In this article we study a diophantine property of probability measures on $\R$. We will always assume that the considered measures have an exponential moment and a drift. We link this property to the points in $\C$ close to the imaginary axis where the Fourier-Laplace transform of those measures take the value $1$ and finally, we apply this to the study of the speed in Kesten's renewal theorem on $\R$.
\end{abstract}

\maketitle

\tableofcontents

\section{Introduction and main results}

In this article, we consider a random walk on $\R$ driven by a probability measure $\rho$  having an exponential moment an a positive drift $\lambda = \int_\R y\di \rho(y)>0$.

The Markov operator associated to $\rho$ is defined, for any borelian non negative function $f$ on $\R$ and any $x\in \R$ by
\[
Pf(x) = \int_\R f(x+y) \di \rho(y)
\]
The question we want to study is : given a function $f \in \cal C^\infty_c(\R)$, can we write $f=g-Pg$ where $g$ is a function on $\R$ vanishing quickly at infinity.

If we set $g=\sum_{n=0}^{+\infty} P^n f$, then $g$ is defined on $\R$ (according to the large deviations inequality) and $g-Pg=f$. Moreover, Kesten's renewal theorem (see~\cite{Fel71}) proves that
\[
\lim_{x\to +\infty} g(x) =0 \text{ and }\lim_{x\to -\infty} g(x) =\frac 1 \lambda \int_\R f(x)\di x
\]
where $\lambda = \int_\R y\di\rho(y)$.

Therefore our question is equivalent to finding the speed in Kesten's theorem.

\medskip
We define the Green kernel $G$ by
\[
G = \sum_{n=0}^{+\infty} \tilde\rho^{\star n} \text{ where for any borelian }A\subset \R,\; \tilde\rho(A) = \int_\R \un_{A}(-y) \di \rho(y)
\]
For $x\in \R$ note 
\[
H(x) = \sum_{n=0}^{+\infty} \rho^{\star n}(]-\infty,x])
\]
Then, we have that for any $f\in \cal S(\R)$,
\[
G\star f(x)= \sum_{n=0}^{+\infty} P^nf(x) = \int_\R f(x+s) H( \di s) = -\int_\R f'(x+s) H(s) \di s
\]
moreover, Smith's renewal theorem (see~\cite{Smi64}) proves that
\begin{equation}\label{equation:renewal_theorem_smith}
\lim_{x\to \pm\infty} H(x) - \left( \frac x \lambda + \frac {\lambda_2}{2\lambda^2} \right) \un_{\R_+}(x) = 0 \text{ where }\lambda_2 = \int_\R y^2 \di\rho(y)
\end{equation}
Thus, we introduce the function $R(x) = H(x) - \left( \frac x \lambda - \frac {\lambda_2}{2\lambda^2} \right) \un_{\R_+}(x)$ and  we get, after some computations, that for any $f\in \cal S(\R)$ and any $x\in \R$,
\begin{equation}\label{equation:control_stieltjes}
\sum_{n=0}^{+\infty} P^nf(x) = \frac 1 \lambda \int_x^{+\infty} f(t)\di t + \frac{\lambda_2}{2\lambda^2} f(x) - \int_\R f'(x+s) R(s) \di s
\end{equation}

In~\cite{Car83}, Carlsson managed to have a control on $R$ assuming some conditions on $\rho$ (a polynomial moment one that holds here since we assume some exponential moment and one he calls ``non lattice of type p'' and which is exactly what we call $p-$weakly diophantine, following Breuillard ans his diophantines measures defined in~\cite{Bre05}).
Using Carlsson's estimates, and equation~\ref{equation:control_stieltjes} we have that the convergence in the renewal theorem is polynomial.

In~\cite{BG07}, Blanchet and Glynn got an exponential control on $R$ but their assumption is to strong for our study since they assume that there is $\varepsilon\in \R_+^\star$ such that
\begin{equation}\label{equation:strongly_non_lattice}
\inf_{|t|\geqslant \varepsilon} |1-\widehat \rho(t)|>0
\end{equation}
This condition is what we will call $0-$weakly-diophantine and our method will prove proposition~\ref{proposition:renewal_theorem_exponential} which says that under this condition, for any $f\in \cal C^2(\R)$ such that $f,f'$ and $f''$ vanish exponential fast at infinity, the speed in the renewal theorem is also exponential. Using their control and equation~\ref{equation:control_stieltjes}, we would have the same result for functions $f\in\cal C^1$ such that $f$ and $f'$ vanishes exponentially fast.

\medskip
In this article, we study an intermediate cases : we want to know if the speed can be faster that what is given by Carlsson's estimates if we take a weaker condition than the strongly non lattice one.

To study the speed in the renewal theorem for measures that are not $0-$weakly-diophantine, we will use the theory of Fourier multipliers on $\R$ and study precisely the function $1/(1-\widehat\rho)$ where $\widehat \rho$ is the Fourier-Laplace transform of $\rho$. We try to find spaces of functions on which the multiplication by $1/(1-\widehat \rho)$ defines a continuous operator.

\medskip
The diophantine condition we make on measures is given in the following
\begin{definition*}[Weakly-diophantine measures]
Let $\rho$ be a (borelian) probability measure on $\R$ and $l\in \R_+$.

We say that $\rho$ is $l-$weakly-diophantine if
\[
\liminf_{b\to \pm\infty} |b|^l|1-\widehat\rho(ib)|>0
\] 
More generally, we say that $\rho$ is weakly-diophantine if it is $l-$weakly-diophantine for some $l\in \R_+$.
\end{definition*}
First, we will prove an alternative definition of weakly-diophantine measures using the zeros of $1-\widehat \rho$. 
\begin{proposition*}[\ref{proposition:zeroes_fourier_transform}]
Let $\rho$ be a probability measure on $\R$ which have an exponential moment of order $\eta$ and a positive drift $\lambda$.

Then, there is $s_0\in \R_+^\star$ such that for any $l\in \R_+$, $\rho$ is $l-$weakly diophantine if and only if there is $C\in \R_+^\star$ such that $0$ is the only zero of $1-\widehat\rho$ in
\[
\left\{z\in \C \middle|\Re(z)\in]-s_0,s_0[\text{ and }\Re(z)\geqslant\frac {-C} {1+|\Im(z)|^l}\right\}
\]
\end{proposition*}

Noting $T_\lambda$ the measure having density $\frac 1 \lambda \un_{\R_-}$ against Lebesgue's measure we have, for any $f\in \cal S(\R)$,
\[
T_\lambda \star f(x) = <T_\lambda \star \delta_x ,f > =\frac{1} \lambda \int_{\R_-} f(x-u) \di u= \frac{1} \lambda \int_{x}^{+\infty} f(u)\di u
\]

We will prove the following characterization of weakly-diophantine measures
\begin{proposition*}[\ref{proposition:tempered_renewal_kernel}]
Let $\rho$ be a probability measure on $\R$ which have an exponential moment and a positive drift $\lambda$.

Then, the two following assertions are equivalent
\begin{enumerate}
\item The measure $\rho$ is weakly diophantine
\item For any $f\in \cal S(\R)$, $(G-T_\lambda)\star f\in\cal S(\R)$
\end{enumerate}
\end{proposition*}

In particular, under the diophantine condition, for any $f\in \cal S(\R)$, the speed in the renewal theorem is faster than any polynomial.

\medskip
As we are interested in faster convergence, we will also prove that if the measure is not $0-$diophantine, then the convergence cannot be exponential :
\begin{proposition*}[\ref{proposition:obstruction_exponential}]
Let $\rho$ be a measure on $\R$ that have an exponential moment and a positive drift $\lambda = \int_\R y\di\rho(y)>0$.

Assume that there is $\gamma\in\R_+^\star$ such that for any $f\in \cal C^\infty_c(\R)$
\[
(G-T_\lambda)\star f(x) \in \cal O\left( e^{-\gamma |x|}\right)
\]
Then, $\rho$ is $0-$weakly-diophantine.
\end{proposition*}

After that, we will take interest in intermediate speed of convergence but we will have to restrict our study to some class of weights function. We will have to control the one-sided Laplace transform of those weight functions, this is why we define,
\begin{equation} \label{eqn:Omega}
\Omega = \left\{\begin{array}{c|l} \multirow{3}{*}{$\omega \in \cal C^0(\R, [1,+\infty[) $}& \omega \text{ is even} \\ &\forall \varepsilon\in \R_+^\star\; \lim_{x\to \pm\infty} e^{-\varepsilon|x|} \omega(x) = 0 \\ &\forall \delta \in \R_+^\star\; \sup_{\substack{z\in \C\\ \Re(z)>0 \text{ and }|z| \geqslant \delta}} \left| \int_{0}^{+\infty} \omega(x) e^{-zx} \di x \right| \text{ is finite}\end{array} \right\}
\end{equation}
This set contains, as we prove in appendix~\ref{appendix:Omega}, $\left(x\mapsto e^{a|x|^\alpha}\right)$ for any $a\in \R_+^\star$ and $\alpha\in ]0,1[$, $\left(x \mapsto  \exp\left(A\ln(1+|x|) (l^{\circ m}(|x|))^M\right) \right)$ for any $A,M\in \R_+$, $m\in \N$ where $l(x) = \ln(1+x)$ and these are good candidates for speeds of convergence between polynomial and exponential.

What we will prove is that either $\omega\in \Omega$ doesn't grow faster than any polynomial or there are functions in $\cal C^\infty_c(\R)$ such that the speed of convergence in the renewal theorem is slower than $\omega$. More precisely we will prove next
\begin{proposition*}[\ref{proposition:not_intermediate_speed}]
Let $\rho$ be a probability measure on $\R$ which have an exponential moment and a negative drift $\lambda$.

Assume that $\rho$ is weakly-diophantine but that there is $l\in \R_+^\star$ such that $\rho$ is \underline{not} $l-$weakly-diophantine.

Let $\omega\in \Omega$ be such that $G-T_\lambda: \left\{\begin{array}{ccc} C^\infty_c(\R) & \to & \mathrm{L}^2_\omega(\R) \\ f & \mapsto & (G-T_\lambda)\star f\end{array}\right.$ is continuous, then $\omega$ doesn't grow faster than any polynomial\footnote{ it means that there is $l\in \N$ such that $\liminf_{x\to +\infty} \frac{\omega(x)}{x^l}=0$.}.
\end{proposition*}
This last proposition means that the speed is not faster than ``faster than any polynomial'' in general for functions in $\cal C^\infty_c(\R)$ if the measure is not $l-$weakly-diophantine for all $l\in \R_+^\star$. And to get faster speeds, we would have to study functions in Gelfand and Shilov's spaces (see~\cite{GC64}) but we won't do so in this article.

\subsection{Notations and assumptions}~

Every (considered) measure on any topological space is borelian.

For $\eta \in \R_+^\star$, we put $\C_\eta = \{z\in \C| |\Re(z)|<\eta \}$ and $\overline{\C_\eta} =  \{z\in \C| |\Re(z)|\leqslant\eta \}$.

For $A,B\subset \R$, we note $A\oplus iB = \{a+ib| a\in A,\; b\in B\}$ in particular, if $A\subset \R$, then $A \oplus i\R = \{z\in \C| \Re(z) \in A\}$

For $f\in \mathrm{L}^1(\R)$, we note $\widehat f$ it's Laplace transform defined for any $z\in \C$ such that the integral is absolutely convergent by
\[
\widehat f(z) = \int_\R f(x) e^{-zx} \di x
\]
then if $\widehat f\in \mathrm{ L}^1(i\R)$ the inversion formula becomes, for a.e. $x\in \R$,
\[
f(x) = \frac 1 {2\pi} \int_\R \widehat f(i\xi) e^{i\xi x} \di \xi
\]
In the same way, for any borelian complex measure $\mu$ of finite total variation and any $z\in \C$ such that the integral is absolutely convergent, we set
\[
\widehat \mu(z) = \int_\R e^{-zx} \di \mu(x)
\]
then, $\widehat \mu$ is also the Fourier-Laplace transform of $\mu$ in the distribution sense and we have that for any $x\in \R$ and any $f\in \cal S(\R)$,
\[
\mu\star f(x) = \int_\R f(x-y) \di \mu(y) = \frac 1 {2\pi} \int_\R e^{i\xi x} \widehat f(i\xi) \widehat \mu(i\xi) \di \xi 
\]

\section{Diophantine properties of measures and their Fourier-Laplace transforms}\label{section:diophantine_measure}

In this section, we study the link between the diophantine properties of a measure (that we will define in subsection~\ref{subsection:diophantine_properties_measures}) and of the set of points close to $i\R$ where it's Fourier-Laplace transform takes the value~$1$.

\medskip
Assume for a moment that $\rho$ is a probability measure on $\R$ that have an exponential moment and a density with respect to Lebesgue's measure. Then, Riemann-Lebesgue's lemma proves that the Fourier-Laplace transform of $\rho$ vanishes at infinity. Thus, since it is uniformly Lipschitzian on some strip containing $i\R$ (as we will see in lemma~\ref{lemma:properties_Fourier_Laplace_transform}), we get that if $\widehat \rho$ doesn't take the value $1$ on $i\R$ except at $0$, then there actually is a strip containing $i\R$ where $\widehat\rho$ takes the value $1$ only at $0$.

In this section, we give another condition, far weaker than having a density with respect to Lebesgue's measure on $\R$ and that gives this sort of control of the points where $\widehat \rho$ takes the value $1$.

\medskip
More specifically, the aim is to state proposition~\ref{proposition:zeroes_fourier_transform} which proves that $\rho$ satisfies to the diophantine condition we define in subsection~\ref{subsection:diophantine_properties_measures} if and only if the Fourier-Laplace transform of $\rho$ doesn't take the value $1$ in a zone of controlled shape except at $0$.

\subsection{The Fourier-Laplace transform of probability measures having an exponential moment}

We start with a few reminds on the Fourier-Laplace transform of probability measures having exponential moments. We give here some properties that are not usually stated and for a more general overview, the reader can refer to the book of Rudin~\cite{Rud91}.

\begin{definition}[Fourier-Laplace transform]
Let $\rho$ be a borelian probability measure on $\R$ which have exponential moments of order at least $\eta$ for some $\eta\in\R_+^\star$.

For $z\in \C_\eta$, we set
\[
\widehat \rho(z) = \int_\R e^{-zy}\di\rho(y)
\]
\end{definition}

As we will see, we have all kinds of control on $\widehat \rho$ and it's derivatives close to the imaginary axis. We sum-up these properties in next

\begin{lemma}\label{lemma:properties_Fourier_Laplace_transform}
Let $\rho$ be a probability measure on $\R$ which have an exponential moment of order at least $\eta$.

Then, $\widehat\rho$ is holomorphic on $\C_\eta$ and for any $k\in \N$, and any $\alpha\in[0,\eta[$, $\widehat\rho^{(k)}$ is bounded and uniformly continuous on $\overline{\C_\alpha}$.

Moreover, for any $z\in \C_\eta$, $\widehat\rho(\overline{z}) = \overline{\widehat\rho(z)}$. In particular, if $z$ is a zero of $1-\widehat\rho$, then so is $\overline z$.

Finally, for any $n\in \N^\star$ and any $z\in \C_\eta$,
\[
\widehat{\rho^{\star n}}(z) = \widehat \rho(z) ^n
\]
\end{lemma}

\begin{proof}
Let $z\in \C_\eta$, then for any $x\in \R$, $|e^{-zx}|= e^{-\Re(z)x}\leqslant e^{\eta |x|}$ and so, as $\rho$ has an exponential moment of order at least $\eta$, $\widehat{\rho}$ is holomorphic on $\C_\eta$. As $\rho$ is a (real) measure on $\R$, it is clear that for any $z\in \C_\eta$ $\widehat\rho(\overline z) = \overline{\widehat\rho(z)}$.

By the theorem of differentiation under the integral sign, we have that for any $k\in \N$, 
\[
\widehat\rho^{(k)}(z) = \int_\R (-x)^k e^{-zx} \di\rho(x)
\]
so, for any $\alpha\in[0,\eta[$ and any $z\in\overline{\C_\alpha}$,
\[
|\widehat\rho^{(k)}(z) | \leqslant \int_\R |x|^k e^{\alpha |x|}\di\rho(x)
\]
which is finite since $\rho$ has an exponential moment of order $\eta$.

This proves that $\rho$ and all it's derivatives are uniformly bounded on $\overline{\C_\alpha}$ and the mean value inequality proves that $\rho$ and it's derivatives are even uniformly Lipschitzian on $\overline{\C_\alpha}$.

Finally, if $z\in \C_\eta$ and $n\in \N$,
\[
\widehat{\rho^{\star n}} (z)=\int_\R e^{-zy} \di\rho^{\star n} (y) = \int_{\R^n} e^{-z(y_1+\dots +y_n)} \di\rho(y_1)\dots\di\rho(y_n) = \widehat\rho(z)^n
\]
(to be more specific, we take first $z\in \{-\eta,\eta\}$ to have a non negative function, use Fubini's theorem in this case and see that $\rho^{\star n}$ also have an exponential moment of order at least $\eta$ and then we do the previous computation for any $z$ in  $\C_\eta$ using Fubini's theorem for absolutely integrable functions).
\end{proof}

\subsection{Diophantine properties of measures on \texorpdfstring{$\R$}{R}}\label{subsection:diophantine_properties_measures}

\subsubsection{Preliminary lemmas}We will see in the sequel that the key ingredient to study the speed in Kesten's renewal theorem is to localize the points $z$ that are close to the imaginary axis and where $\widehat\rho(z)=1$ (this set is discrete since $\widehat\rho$ is analytic on $\C_\eta$ as we saw in lemma~\ref{lemma:properties_Fourier_Laplace_transform}). In this part, we will see two important results. On the one hand, lemma~\ref{lemma:approximate_zeros} together with remark~\ref{remark:approximate_zeros}, shows that there always is a solution of $\widehat\rho(z)=1$ close to an approximate one. On the other hand, in lemma~\ref{lemma:almost_lattice} we will see that those points $z$ such that $\widehat\rho(z)=1$ are almost a lattice in $\R^2$ in some sense that will be made precise.

\begin{definition}
Let $\rho$ be a probability measure on $\R$ which has a finite first moment.

We say that $\rho$ has a drift if $\int_\R y\di\rho(y)\not=0$.
\end{definition}

\begin{remark}
In this article we will focus on measures having a drift. We will always assume it to be positive, but, of course, if it is negative, then the measure $\tilde \rho$ defined by $\tilde\rho(A) = \int_\R \un_A(-y) \di\rho(y)$ has a positive one and it's Fourier-Laplace transform is the image of the one of $\rho$ by the application $z\mapsto -z$ on $\C$. Therefore, anything done for measures with positive drifts can be translated for measures with negative ones.
\end{remark}

\begin{lemma}\label{lemma:approximate_zeros}
Let $\rho$ be a probability measure on $\R$ which have an exponential moment of order $\eta$ and a drift $\lambda\not=0$.

For $z\in \C_{\eta/2}$, note
\[
\rho^\circ (z) = \left(\int_\R \left| e^{-zy}-1 \right|^2 \di\rho(y) \right)^{1/2}
\]

Then, for any $\eta'<\eta$, there are $\varepsilon_0,C\in \R_+^\star$ such that for any $z_0\in \C_{\eta'/2}$, if
\[
\rho^{\circ }(z_0)\leqslant \varepsilon_0
\]
there is $a\in \C_{\eta'}$ such that $\widehat \rho(a)=1$ and
\[
|z_0-a| \leqslant C \rho^{\circ}(z_0)
\]
\end{lemma}

\begin{remark}\label{remark:approximate_zeros}
If $z=ib\in i\R$, then
\[
\rho^\circ(ib)^2 = \int_\R \left|e^{-iby}-1\right|^2 \di\rho(y) = 2-2\Re(\widehat \rho(ib)) \leqslant 2|1-\widehat \rho(ib)|
\]
and the previous lemma proves that if $b$ is such that $|1-\widehat \rho(ib)|$ is small, then there is $z$ close to $ib$ and such that $\widehat\rho(z)=1$.
\end{remark}

\begin{proof}
If $\rho^\circ(z_0)=0$, then for $\rho-$a.e. $y\in \R$, $e^{-zy}=1$ and so, $\widehat \rho(z) = 1$ and the lemma is proved with $a=z_0$.
So we may assume without any loss of generality that $\rho^\circ(z_0) \not=0$.

Let $\eta'\in ]0,\eta[$. For $z\in \C_{\eta'/2}$, note $f(z) = 1-\widehat \rho(z)$.
As $0$ is a simple zero of $f$, the argument principle shows that for $r\in\R_+^\star$ small enough,
\[
1 = \frac 1 {2i\pi} \int_{|z|=r} \frac{f'(z)}{f(z)} \di z
\]
Moreover, we also know that if $z\mapsto f(z+z_0)$ doesn't vanish on the circle of radius $r$, then the number of it's zeros in the disk is also
\[
\frac 1 {2i\pi}\int_{|z|=r} \frac{f'(z+z_0)}{f(z+z_0)} \di z
\]
Thus, to prove the lemma, we are going to show that under it's assumptions, $f(z+z_0)$ doesn't vanish on the circle of radius $r$ and that 
\[
\frac 1 {2\pi} \left|\int_{|z|=r} \frac{f'(z)}{f(z)} - \frac{f'(z+z_0)}{f(z+z_0)} \di z \right| <1
\]
But, for $z$ such that $f(z+z_0)\not=0$, we have that
\[
 \frac{f'(z)}{f(z)} - \frac{f'(z+z_0)}{f(z+z_0)} = \frac{f'(z)f(z+z_0)- f'(z+z_0) f(z)}{f(z) f(z+z_0)}
\]
And,
\begin{equation}\label{equation:approximate_zero_minor}
|f(z+z_0)| \geqslant |f(z)| - |f(z+z_0)-f(z)|
\end{equation}
Moreover,
\begin{flalign*}
|f(z+z_0)-f(z)| &= \left| \int_\R e^{-zy} (e^{-z_0y}-1) \di \rho(y)\right| \\&\leqslant \left(\int_\R e^{-2\Re(z) y}\di\rho(y) \right)^{1/2} \left(\int_\R \left| e^{-z_0y}-1\right|^2 \di \rho(y) \right)^{1/2} \\
& \leqslant C_0\rho^{\circ}(z_0) \text{ where we put }C_0=\left(\int_\R e^{\eta'|y|} \di \rho(y)\right)^{1/2}
\end{flalign*}
And, as $f(0)=0$, we have that
\[
\frac{f(z)}{z} \xrightarrow[z\to 0]\, f'(0) = \int_\R y \di\rho(y) \not=0
\]
Thus, there is $\varepsilon_1$ such that for any $z\in B(0,\varepsilon_1)$,
\[
|f(z)| \geqslant |z| \frac{|\lambda|}{2}
\]
And so, for any $z\in B(0,\varepsilon_1)$, equation~\ref{equation:approximate_zero_minor} becomes
\[
|f(z)||f(z+z_0)| \geqslant |z| \frac{|\lambda|}{2} \left(|z| \frac{|\lambda|}{2} -C_0\rho^\circ(z_0) \right)
\]

Then,
\begin{flalign*}
f'(z) f(z+&z_0) - f'(z+z_0)f(z) = \widehat \rho'(z) - \widehat \rho'(z+z_0) - \widehat \rho'(z) \widehat\rho(z+z_0) + \widehat\rho'(z+z_0) \widehat \rho(z) \\
&= \int_\R ye^{-zy} \left( e^{-z_0y}-1\right)\di \rho(y) + \int_{\R^2} y_1 e^{-z(y_1+y_2)} \left(e^{-z_0y_1} - e^{-z_0y_2} \right) \di \rho^{\otimes 2}(y_1,y_2)
\end{flalign*}
but,
\begin{flalign*}
\left| \int_\R ye^{-zy} \left( e^{-z_0y}-1\right)\di \rho(y)\right| &\leqslant \left(\int_\R y^2e^{-2\Re(z)y}\di\rho(y)\right)^{1/2}  \left(\int_\R \left| e^{-z_0y}-1\right|^2 \di\rho(y) \right)^{1/2} \\
& \leqslant C_1\rho^\circ(z_0) \text{ where we put }C_1=\left(\int_\R y^2 e^{\eta'|y|} \di \rho(y)\right)^{1/2}
\end{flalign*}
And,
\begin{flalign*}
\bigg| \int_{\R^2} y_1  e^{-z(y_1+y_2)} & \left(e^{-z_0y_1} - e^{-z_0y_2} \right) \di \rho^{\otimes 2}(y_1,y_2) \bigg|^2 \\
& \leqslant \int_{\R^2} y_1^2 e^{-2\Re(z)(y_1+y_2)} \di \rho^{\otimes 2}(y_1,y_2)  \int_{\R^2} \left|e^{-z_0y_1} - e^{-z_0y_2} \right|^2 \di \rho^{\otimes 2}(y_1,y_2) \\
& \leqslant C_0^2 C_1^2 \int_{\R^2} \left|e^{-z_0y_1} - e^{-z_0y_2} \right|^2 \di \rho^{\otimes 2}(y_1,y_2)
\end{flalign*}
but,
\begin{flalign*}
\left(\int_{\R^2} \left|e^{-z_0y_1} - e^{-z_0y_2} \right|^2  \di \rho^{\otimes 2}(y_1,y_2)\right)^{1/2} &\leqslant \left(\int_{\R^2} \left|e^{-z_0y_1} - 1 \right|^2 \di \rho^{\otimes 2}(y_1,y_2) \right)^{1/2} \\
& \retrait +\left(\int_{\R^2} \left|e^{-z_0y_2} - 1 \right|^2  \di \rho^{\otimes 2}(y_1,y_2)\right)^{1/2}   \\
&\leqslant 2 \rho^\circ(z_0)
\end{flalign*}

What we just prove is that for any $z$ such that $|z|=r<\varepsilon_1$ we have that
\[
\left| \frac{f'(z)}{f(z)} - \frac{f'(z+z_0)}{f(z+z_0)}  \right| \leqslant \frac{C_1(C_0+1) \rho^{\circ}(z_0)}{r(r-C_0\rho^{\circ}(z_0))}
\]
So, if $r=C_2\rho^{\circ}(z_0) $ with $C_2>C_0$ such that $C_1(C_0+1)/(C_2-C_0)<1$, and $\frac{\varepsilon_1}{C_2}= \varepsilon_0$, we get that
\[
\left|\frac 1 {2i\pi}\int_{|z|=r} \frac{f'(z)}{f(z)} - \frac{f'(z+z_0)}{f(z+z_0)}  \di z\right| \leqslant \frac{C_1(C_0+1) \rho^{\circ}(z_0)}{r-C_0\rho^{\circ}(z_0)} = \frac{C_1(C_0+1) }{C_2-C_0}<1
\]
and this finishes the proof of the lemma.
\end{proof}

In next lemma, we prove that there is almost an additive structure on the points $z$ where $\rho(z)=1$. In particular, it proves that there are plenty of them.

\begin{lemma}\label{lemma:almost_lattice}
Let $\rho$ be a probability measure on $\R$ which have an exponential moment of order $\eta$ and a negative drift $\lambda$.

There are $\varepsilon_0\in \R_+^\star$ and $C\in \R$ such that if $z_1,z_2\in \C_{\varepsilon_0}$ are such that $\widehat\rho(z_1)=1=\widehat\rho(z_2)$. Then, there is $z_3\in \C_\eta$ such that $\widehat \rho(z_3) =1 $ and
\[
|z_3-(z_1+z_2)| \leqslant C \left( |\Re(z_1)|+ \sqrt{|\Re(z_2)|} \right)^{1/2}
\]
\end{lemma}

\begin{proof}

\[
\widehat\rho(z_1+z_2) - 1 = \int_\R e^{-(z_1+z_2)y}  \di\rho(y) - \int_\R e^{-z_1 y} \di \rho(y) = \int_\R e^{-z_1y} \left( e^{-z_2y} -1 \right) \di\rho(y)
\]
So,
\begin{flalign*}
\left| \widehat\rho(z_1+z_2) - 1\right| &\leqslant \int_\R e^{-2 \Re(z_1)y} \left| e^{-z_2y} -1 \right| \di\rho(y)\\
&\leqslant \left( \int_\R e^{-2\Re(z_1)y} \di\rho(y)  \right)^{1/2} \left(\int_\R  \left| e^{-z_2y} -1 \right|^2 \di\rho(y)\right)^{1/2}
\end{flalign*}
But, $ \widehat\rho(z_2)=1$, so
\[
\int_\R  \left| e^{-z_2y} -1 \right|^2 \di\rho(y) = 1 + \int_\R e^{-2\Re(z_2)y} \di \rho(y) -2 \Re(\widehat \rho(z_2)) = \int_\R e^{-2\Re(z_2)y} \di \rho(y)-1
\]
So,
\begin{flalign*}
\rho^\circ(z_1&+z_2)^2 = \int_\R \left| e^{-(z_1+z_2)y}-1 \right|^2 \di\rho(y) = \int_\R e^{-2\Re(z_1+z_2)y}\di\rho(y) +1 -2\Re(\widehat\rho(z_1+z_2)) \\
&\leqslant \left|\int_\R e^{-2\Re(z_1+z_2)y}\di\rho(y) -1 \right| + 2 \left( \int_\R e^{-2\Re(z_1)y} \di\rho(y)  \right)^{1/2} \left| \int_\R e^{-2\Re(z_2)y} \di \rho(y)-1\right|^{1/2}
\end{flalign*}

Finally, lemma~\ref{lemma:properties_Fourier_Laplace_transform} proves that $\rho$ is uniformly Lipschitzian on $\C_{\varepsilon_0}$ and so, there is a constant $C_0$ such that
\[
\left|\int_\R e^{-2\Re(z_1+z_2)y}\di\rho(y) -1 \right| \leqslant C_0\left( |\Re(z_1)| + |\Re(z_2)| \right)
\]
and
\[
2 \left( \int_\R e^{-2\Re(z_1)y} \di\rho(y)  \right)^{1/2} \left| \int_\R e^{-2\Re(z_2)y} \di \rho(y)-1\right|^{1/2} \leqslant C_0\sqrt{|\Re(z_2)|} 
\]
This proves that
\[
\rho^\circ(z_1+z_2) \leqslant 2C_0 (|\Re(z_1)| + \sqrt{|\Re(z_2)|})
\]
So, if $\Re(z_1)$ and $\Re(z_2)$ are small enough, we can use lemma~\ref{lemma:approximate_zeros} to conclude.
\end{proof}

\begin{remark}\label{remark:Delone_set}
The previous lemma also shows that if we put $\cal A = \{z\in \C_\eta |\widehat \rho(z)=1\}$ then
\[
\sup_{a\in \cal A} \di(a,\cal A \setminus\{a\}) \text{ is finite}
\]
and this, with lemma~\ref{lemma:properties_Fourier_Laplace_transform}, proves that $\cal A$ is a Delone set (a set that is both relatively dense and uniformly discrete).
\end{remark}

\begin{corollary}
Let $\rho$ be a probability measure on $\R$.

Note
\[
\widehat G_\rho = \left\{ b\in \R\; \widehat \rho(ib)=1\right\} = \left\{b\in \R \middle| \text{ for }\rho-\text{a.e. }y\in \R \;by\in 2\pi \Z \right\}
\]
Then, $\widehat G_\rho$ is a closed discrete subgroup of $\R$.
\end{corollary}

\begin{proof}
The fact that $\left\{ b\in \R\; \widehat \rho(ib)=1\right\}$ is a group can be seen has a straightforward corollary of lemma~\ref{lemma:almost_lattice} if $\rho$ has an exponential moment and a negative drift but it doesn't actually uses these assumptions.

It's easy to see that $\left\{b\in \R \middle| \text{ for }\rho-\text{a.e. }y\in \R \;by\in 2\pi \Z \right\}$ is a group so we will prove that if is equal to $\left\{ b\in \R\; \widehat \rho(ib)=1\right\}$.

Indeed, if $b\in \R$ is such that $\widehat\rho(ib) = 1$, then
\[
\int_\R \left| e^{-iby} -1 \right|^2 \di\rho(y) = 2- 2\Re\widehat\rho(ib) =0
\]
So, for $\rho-$a.e $y\in \R$, $by\in 2\pi \Z$ and this finishes the proof since the other inclusion is easy to get.
\end{proof}

\begin{definition}[Lattice measures]\label{definition:lattice_measure}
Let $\rho$ be a probability measure on $\R$.

We say that $\rho$ is lattice if $\widehat G_\rho \not=\{0\}$.
\end{definition}

As we already saw, a probability measure $\rho$ is lattice if and only if there is $b\in \R^\star$ such that $\rho(\frac{2\pi} b \Z) =1$ and this is the usual definition of lattice probability measures.

\subsubsection{Diophantine measures}~

In this article, we won't do more with lattice measures. Therefore, we will only work with measures such that for any $b\in \R^\star$, $\widehat\rho(ib)\not=1$.

We say that the measure $\rho$ satisfies to Cram\'er's condition or is strongly non lattice if
\[
\liminf_{b\to +\infty} |1-\widehat \rho(ib)|>0
\]
Another corollary of lemma~\ref{lemma:approximate_zeros} is that if $\rho$ satisfies to Cram\'er's condition, then there is a strip containing $i\R$ in $\C$ on which $\widehat\rho-1$ doesn't vanish except at $0$.

This condition is satisfied, as an example, by measures having a density against Lebesgue's measure (according to Riemann-Lebesgue's lemma) or more generally by spread-out measures. It has been studied by Blanchet and Glynn in~\cite{BG07} but we want to relax it to be able to deal with measures such as $\frac 1 2(\delta_1 + \delta_a)$ for some diophantine number $a$.

Thus, we introduce a class of measures such that we may have
\[
\liminf_{b\to +\infty} |1-\widehat \rho(ib)|=0
\]
but the convergence is not too quick.

Breuillard defined in~\cite{Bre05} the diophantine measures. These are measures such that for some $l\in \R_+$,
\[
\liminf_{b\to \pm\infty} |b|^l (1-|\widehat\rho(ib)|)>0
\]
This condition is once again too restrictive for our study because we really don't need to control $1-|\widehat{\rho}|$ but just $1-\widehat\rho$ and so, we make the following
\begin{definition}[Weakly-diophantine measures]
Let $\rho$ be a (borelian) probability measure on $\R$ and $l\in \R_+$.

We say that $\rho$ is $l-$weakly-diophantine if
\[
\liminf_{b\to \pm\infty} |b|^l|1-\widehat\rho(ib)|>0
\] 
More generally, we say that $\rho$ is weakly-diophantine if it is $l-$weakly-diophantine for some $l\in \R_+$.
\end{definition}

\begin{remark}
If $\rho$ is weakly-diophantine, then it isn't lattice since there can't be infinitely many $b\in \R$ such that $\widehat \rho(ib)=1$.
\end{remark}

The following lemma corresponds to the proposition~3.1 of~\cite{Bre05}. 
\begin{lemma}\label{lemma:equivalent_definition_diophantine_measure}
Let $\rho$ be a probability measure on $\R$ and $l\in \R_+$.

Define the three assertions
\begin{enumerate}[label=(\roman*)]
\item \label{item:definition} The measure $\rho$ is $l-$weakly-diophantine
\item \label{item:distance_integers} \[
\liminf_{b\to \pm\infty} |b|^l \int_\R \{b y\}^2 \di\rho(y) >0 \text{ where for }t\in \R,\; \{t\} = \inf_{p\in \Z} |t-p| = d(t,\Z)>0
\]
\item \label{item:lattice}
\[
\liminf_{b\to \pm\infty} |b|^l \int_\R \left|e^{-iby} -1 \right|^2 \di\rho(y) > 0
\]
\end{enumerate}
Then, \ref{item:distance_integers} is equivalent to \ref{item:lattice} with the same $l$.

Moreover, if \ref{item:distance_integers} or equivalently \ref{item:lattice} hold then we have \ref{item:definition} with the same $l$. 

Finally, if $\rho$ is $l-$diophantine, then \ref{item:distance_integers} or equivalently \ref{item:lattice} hold with $2l$.
\end{lemma}

\begin{proof}
First of all, remark that for any $b\in \R$,
\[
\frac 1 2 \int_\R \left|e^{-iby}-1 \right|^2 \di \rho(y) = 1-\Re(\widehat \rho(ib)) = \int_\R 1-\cos(by) \di \rho(y)
\]

As Breuillard, we note that there are $c_1,c_2\in \R_+^\star$ such that for any $t\in \R$,
\[
c_1\{t\}^2 \leqslant 1-\cos(2\pi t) \leqslant c_2\{t\}^2 \text{ where we noted }\{t\} = d(t,\Z)
\]
and this proves that \ref{item:distance_integers}$\Leftrightarrow$\ref{item:lattice}.

\medskip
Moreover, for any $z\in \C$ such that $|z|\leqslant 1$,
\[
\frac 1 2|z-1|^2 = \frac 1 2 \left(1 +|z|^2 -2\Re(z) \right)\leqslant \left|1-\Re(z)\right| \leqslant |1-z|
\]
So, for any $b\in \R$,
\[
\frac 1 2|1-\widehat \rho(ib)|^2 \leqslant \left|1-\Re(\widehat\rho(ib)) \right| \leqslant |1-\widehat \rho(ib)|
\]
and this proves the end of the lemma.
\end{proof}

We give some examples of diophantine measures and some ways to construct them in next
\begin{lemma}\label{lemma:examples_diophantine_measures}
Let $\rho$ be a borelian probability measure on $\R$.
\begin{enumerate}[label=(\roman*)]
\item \label{item:spread_out} If $\rho$ is spread-out, then it is $0-$weakly diophantine
\item \label{item:diophantine} If $\rho$ is $l-$diophantine, it is also $l-$weakly-diophantine.
\item \label{item:diophantine_number} If $a$ is a diophantine number, then $\frac 1 2 \delta_1 + \frac 1 2 \delta_a$ is weakly-diophantine but not diophantine. In particular, for a.e. $x,y\in \R$, $\frac 1 2 \delta_x + \frac 1 2 \delta_y$ is weakly-diophantine but not diophantine.
\item \label{item:density}If for some $\delta\in\R_+^\star$ and some $l-$weakly diophantine probability measure $\rho_1$ we have that $\rho-\delta \rho_1$ is a non negative measure, then $\rho$ is also $l-$weakly-diophantine.
\item \label{item:power} If $\rho^{\star m}$ is $l-$weakly-diophantine, then so does $\rho$.
\end{enumerate}
\end{lemma}

\begin{proof}
\begin{enumerate}
\item[\ref{item:spread_out}] This is Riemann-Lebesgue's lemma.
\item[\ref{item:diophantine}]The triangular inequality proves that for any $t\in \R$,
\[
|1-|\widehat \rho(it)|| \leqslant |1-\widehat \rho(it)|
\]
and that proves~\ref{item:diophantine}.
\item[\ref{item:diophantine_number}] This is the definition of being a diophantine number applied with lemma~\ref{lemma:equivalent_definition_diophantine_measure}.
\item[\ref{item:density}] Let $\rho_1$ be a $l-$weakly-diophantine measure and $\delta\in\R_+^\star$ such that $\rho-\delta\rho_1$ is positive.
Then for any $t\in \R$,
\[
\int_\R \left| 1-e^{-ity}\right|^2 \di\rho(y) \geqslant \delta\int_\R \left| 1-e^{-ity}\right|^2 \di\rho_1(y)
\]
and we conclude with lemma~\ref{lemma:equivalent_definition_diophantine_measure}
\item[\ref{item:power}] For any $t\in \R$ and any $m\in \N^\star$, we have that
\[
\left|1-\widehat{\rho^{\star m}}(it)\right| = \left|1-\widehat \rho(it) ^m \right| = \left|1-\widehat \rho(it)\right| \left| 1 + \dots + \widehat \rho(it)^{m-1} \right|  \leqslant m \left|1-\widehat \rho(it)\right|
\]
so if $\rho^{\star m}$ is $l-$weakly-diophantine, so does $\rho$.
\end{enumerate}
\end{proof}

\begin{lemma}\label{lemma:strong_weak_diophantine}
If $\rho$ is symmetric and $l-$weakly-diophantine, then it is also $2l-$strongly-diophantine.
\end{lemma}

\begin{example}
Breuillard gave in~\cite[Example 3.1, 3.2]{Bre05} an example of a measure $\rho_1$ which is not diophantine and such that $\supp(\rho) \supset[0,1]$. He even gave an example $\rho_2$ of such a measure without atoms. The measures $\rho_1 \star \tilde\rho_1$ and $\rho_2 \star \tilde\rho_2$ are not weakly-diophantine according to lemma~\ref{lemma:strong_weak_diophantine} even if $\supp \rho_1 \star \tilde\rho_2$ has a non empty interior and $\rho_2\star\tilde\rho_2$ has no atoms.
\end{example}

\begin{proof}[Proof of lemma~\ref{lemma:strong_weak_diophantine}]
Since $\rho$ is symmetric, we have that for any $t\in\R$, $\widehat \rho(ib) \in \R$. Therefore, $|\widehat \rho(ib)|^2 = \widehat\rho(ib)^2$.

Assume that $\rho$ is not $2l-$diophantine. It means that
\[
\liminf_{b\to \pm\infty} |b|^{2l} |1-|\widehat \rho(ib)|| = 0
\]
And so, we have a sequence $(b_n)$ converging to infinity and such that
\[
\lim_{n\to +\infty} |b_n|^{l} |1-\widehat \rho(ib_n)| = 0 \text{ or }\lim_{n\to +\infty} |b_n|^{l} |1+\widehat \rho(ib_n)| = 0
\]
The first alternative means that $\rho$ is not $2l-$weakly diophantine and, in particular, it isn't $l-$weakly-diophantine and we are going to prove that so does the second alternative.

Indeed, for any $t\in \R$,
\[
2|1+ \widehat \rho(it)| \geqslant 2|1+ \Re(\widehat\rho(it))| = \int_\R \left| e^{-ity}+1 \right|^2 \di \rho(y)
\]
and,
\[
\int_\R \left| e^{-2ity} -1 \right|^2 \di\rho(y) = \int_\R \left| e^{-ity} -1 \right|^2\left| e^{ity} +1 \right|^2 \di\rho(y) \leqslant 2\int_\R \left| e^{-ity} +1 \right|^2 \di\rho(y)
\]
Therefore,
\[
\lim_{n\to +\infty} |b_n|^{2l} |1+\widehat \rho(ib_n)| = 0 \Rightarrow \lim_{n\to +\infty} |2b_n|^{2l} \int_\R \left| e^{-i(2b_n)y} -1 \right|^2 \di\rho(y)=0
\]
and this means, with lemma~\ref{lemma:equivalent_definition_diophantine_measure}, that $\rho$ is not $l-$weakly-diophantine.
\end{proof}

\subsection{Shape of the zero-free region}
In this section, we go back to the study of points where the Fourier Laplace transform of $\rho$ takes the value $1$ when $\rho$ has a positive drift.
First, we study zeros of $1-\widehat\rho$ on a strip on the right of the imaginary axis in the complex plane. We actually show that there are none except $0$. We also prove that the zeroes of $\rho$ on a strip on the left of the imaginary axis are uniformly isolated.

\begin{lemma}\label{lemma:uniformily_isolated_zeros}
Let $\rho$ be a probability measure on $\R$ which have an exponential moment and a positive drift $\lambda=\int_\R y\di\rho(y)>0$.

Then, there is $s_0\in \R_+^\star$ such that for any $z\in ]0,s_0]\oplus i\R $, $|\widehat\rho(z)|<1$.

Moreover, we also have that the zeros of $1-\widehat\rho$ in $\C_{s_0}$ are uniformly isolated i.e. :
\[
\inf_{\substack{z_1,z_2\in \C_{s_0}\\ \widehat\rho(z_1)=1=\widehat\rho(z_2)\\z_1\not=z_2}} |z_1-z_2|>0
\]
\end{lemma}

\begin{proof}
For $s\in ]-\eta,\eta[$ and $b\in \R$, we have that
\[
\left|\widehat \rho(s+ib)\right| = \left| \int_\R e^{-(s+ib) y}\di\rho(y) \right| \leqslant \int_\R e^{-sy} \di\rho(y) = \widehat \rho(s)
\]
Moreover, since $\rho$ has exponential moments, $\widehat \rho$ is defined on a neighbourhood of $0$ and is differentiable at $0$. Then, $\widehat\rho(0)=1$ and $\widehat\rho'(0)=-\int_\R y\di\rho(y) =-\lambda<0$ and so, there is $s_0\in\R_+^\star$ such that for any $s\in ]0,s_0[$, $\widehat\rho(s)<1$.

Thus, for any $z\in]0,s_0[\oplus i\R$, $|\widehat \rho(z)|<1$ (and in particular, $\widehat\rho(z)\not=1$).

\medskip
To prove that the zeros of $1-\widehat\rho$ are uniformly isolated, note $f=\widehat\rho-1$. Then, $f$ is holomorphic on $\C_\eta$ and for any $z\in \C_\eta$,
\[
f'(z)=- \int_\R y e^{-zy}\di\rho(y) \text{ and }|f''(z) |=\left| \int y^2 e^{-zy}\di\rho(y) \right| \leqslant \int_\R y^2e^{-\Re(z) y} \di\rho(y)
\]
And so, $f''$ is uniformly bounded on $\C_\eta$ since $\rho$ has exponential moments. Thus, $f'$ is uniformly Lipschitzian on $\C_\eta$.

Moreover, we noted $\lambda = \int_\R y\di\rho(y)$ so that, for $z\in \C_\eta$,
\begin{flalign*}
|f'(z)+\lambda|^2 &= \left|\int_\R y(e^{-zy}-1)\di\rho(y) \right|^2 \leqslant \int_\R y^2\di\rho(y) \int_\R |e^{-zy} -1|^2 \di\rho(y) \\
&\leqslant \int_\R y^2 \di\rho(y) \left( 1+ \widehat\rho(2\Re(z))-2\Re(\widehat\rho(z)) \right)
\end{flalign*}
So, if $z$ is such that $\widehat\rho(z)=1$, then
\[
|f'(z)+\lambda|^2 \leqslant\int_\R y^2 \di\rho(y) \left( \widehat\rho(2\Re(z))-1 \right) 
\]
and, as $\widehat \rho-1$ is continuous and vanishes at $0$, we may decrease $s_0$ so that for any $z\in \C_{s_0}$ such that $f(z)=0$, $|f'(z)+\lambda|<|\lambda|/2$ and as $f'$ is uniformly continuous according to lemma~\ref{lemma:properties_Fourier_Laplace_transform}, this tells us that there are $\varepsilon,\delta$ such that for any $z\in \C_{s_0}$ such that $f(z)=0$ and any $z'\in \cal B(z,\varepsilon)$, $|f'(z)|>\delta$.

Thus, if $z_0\in \C_{s_0}$ is such that $f(z_0)=0$, then for any $z_1,z_2 \in \cal B(z_0,\varepsilon)$, the mean value inequality says that
\[
|z_1-z_2| \leqslant |f(z_1)-f(z_2)| \sup_{z\in B(z_0,\varepsilon)} \frac 1 {|f'(z)|} \leqslant \frac 1 {\delta} |f(z_1)-f(z_2)|
\]
This means that $f$ vanishes only once in $B(z_0,\varepsilon)$ and this is what we intended to prove.
\end{proof}

We are now ready to study the shape of the region where $\widehat\rho$ doesn't take the value $1$.

\begin{proposition}\label{proposition:zeroes_fourier_transform}
Let $\rho$ be a probability measure on $\R$ which have an exponential moment of order $\eta$ and a positive drift $\lambda$.

Then, there is $s_0\in \R_+^\star$ such that for any $l\in \R_+$, $\rho$ is $l-$weakly diophantine if and only if there is $C\in \R_+^\star$ such that $0$ is the only zero of $1-\widehat\rho$ in
\[
\left\{z\in \C \middle|\Re(z)\in]-s_0,s_0[\text{ and }\Re(z)\geqslant\frac {-C} {1+|\Im(z)|^l}\right\}
\]
\end{proposition}

Before the proof, we draw this zero-free region.
\begin{figure}[H]\label{picture:zero-free_region}
\centering
\begin{tikzpicture}
\fill[pattern = north west lines] (1,-2) -- (1,2) -- plot [domain=2:-2,variable=\t,samples=200] ({max(-1.1/(1+(abs(\t))^(2)),-1)},\t) -- cycle ;
\draw [very thin,->] (-2,0) -- (2,0) ;
\draw [very thin,->] (0,-2.2) -- (0,2.2) ;
\draw (-0.2,0) node[below]{$O$} ;
\draw (1,0.2) node[right]{$s_0$} ;
\draw (-1,0.2) node[left]{$-s_0$} ;
\draw (3,0) node[right]{$\left\{z\in \C \middle|\Re(z)\in]-s_0,s_0[\text{ and }\Re(z)\geqslant\frac {-C} {1+|\Im(z)|^l}\right\}$} ;
\draw [ultra thick] (1,-2) -- (1,2) ;
\draw [ultra thick] plot [domain=-2:2,variable=\t,samples=200] ({max(-1.1/(1+(abs(\t))^(2)),-1)},\t); 
\end{tikzpicture}
\caption{Shape of the zero-free region of $1-\widehat\rho$ for a diophantine measure $\rho$}
\end{figure}
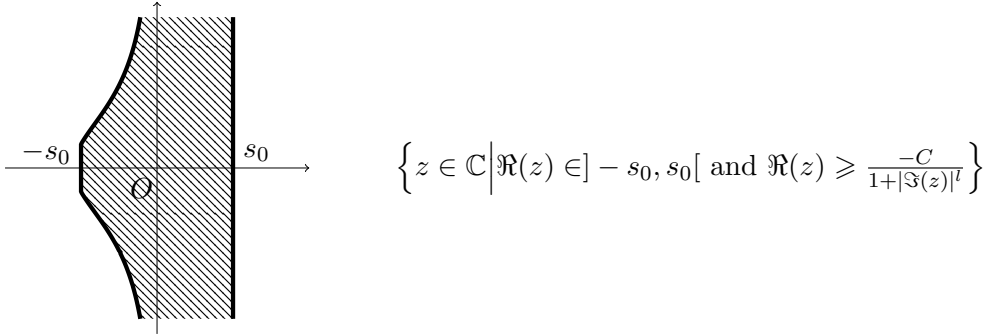

We divide the proof in two lemmas, each one proving one implication.

\begin{lemma}\label{lemma:zeroes_fourier_transform_if_part}
Let $\rho$ be a probability measure on $\R$ which have an exponential moment of order $\eta$ and a positive drift $\lambda$.

Then, there is $s_0\in \R_+^\star$ such that for any $l\in \R_+$, if $\rho$ is $l-$weakly diophantine then there is $C\in \R_+^\star$ such that $0$ is the only zero of $1-\widehat\rho$ in
\[
\left\{z\in \C \middle|\Re(z)\in]-s_0,s_0[\text{ and }\Re(z)\geqslant\frac {-C} {1+|\Im(z)|^l}\right\}
\]
\end{lemma}

\begin{proof}
Assume that there is no such $C$. This means that for any $\varepsilon\in\R_+^\star$, there is $z_\varepsilon\in \C_\eta\setminus\{0\}$ such that
\[
\widehat\rho(z_\varepsilon)=1\text{ and }\Re(z_\varepsilon) \geqslant \frac{ -\varepsilon}{1+|\Im(z_\varepsilon)|^l}
\]
Lemma~\ref{lemma:uniformily_isolated_zeros} proves also that $\Re(z_\varepsilon)\leqslant 0$ and that $\lim_{\varepsilon\to 0^+} \Re(z_\varepsilon)=0$ thus,
\[
\lim_{\varepsilon\to 0^+} |\Im(z_\varepsilon)|=+\infty
\]

Then,
\[
|1- \widehat \rho(i\Im(z_\varepsilon))|= |\widehat\rho(z_\varepsilon) - \widehat \rho(i\Im(z_\varepsilon))| \leqslant |\Re(z_\varepsilon)| \sup_{z\in \C_\eta} |\widehat\rho'(z)| \leqslant \frac{\varepsilon}{1+|\Im(z_\varepsilon)|^l} \sup_{z\in \C_\eta} |\widehat\rho'(z)|
\]
where $\sup_{z\in \C_\eta} |\widehat \rho'(z)|$ is finite according to lemma~\ref{lemma:properties_Fourier_Laplace_transform}.

So,
\[
\limsup_{\varepsilon\to 0^+}|\Im(z_\varepsilon)|^l |1-\widehat\rho(i\Im(z_\varepsilon))| \leqslant \limsup_{\varepsilon\to 0^+} \frac{\varepsilon|\Im(z_\varepsilon)|^l}{1+|\Im(z_\varepsilon)|^l} \sup_{z\in \C_\eta} |\widehat\rho'(z)| =0
\]
And $\rho$ is not $l-$weakly diophantine.
\end{proof}

\begin{lemma}\label{lemma:zeroes_fourier_transform_only_if_part}
Let $\rho$ be a probability measure on $\R$ which have an exponential moment of order $\eta$ and a positive drift $\lambda$.

Then there is $\eta_0$ such that if $\rho$ is not $l-$weakly-diophantine, we have that for any $\varepsilon\in\R_+^\star$, there is $(z_\varepsilon)\in \C_{\eta_0}$ such that
\[
\widehat\rho(z_\varepsilon) =1\text{ and }\Re(z_\varepsilon) \geqslant \frac{-\varepsilon}{1+|\Im(z_\varepsilon)|^l}
\]
\end{lemma}

\begin{remark} \label{remark:invariant_measures}
The previous lemma proves that if $\rho$ is not $0-$weakly diophantine, then there are non constant $\rho-$harmonic function on $\R$ of arbitrarily slow exponential growth.
This means that for any $\varepsilon \in \R_+^\star$, there is a non constant (continuous) function $g$ on $\R$ such that $\rho \star g=g$, $\sup_{x\in \R} e^{-\varepsilon|x|}|g(x)|$ is finite and $\limsup_{x\to -\infty} \frac 1 {|x|} \ln|g(x)| >0$. We just take $g(x) = e^{zx}$ where $z$ is such that $\widehat \rho(z)=1$ and $|\Re(z)| \leqslant \varepsilon$. Or, if we want a real-valued function, we take $g(x) = \Re(e^{zx})$.
\end{remark}

\begin{proof}
We will actually prove that there are in fact infinitely many such $z_\varepsilon$.

Assume that $\rho$ is not $l-$diophantine. This means that there is $(b_n)\in \R^\N$ such that $|b_n| \xrightarrow\, +\infty$ and
\begin{equation}\label{equation:lemma_zeros}
|b_n|^l|1-\widehat\rho(ib_n)| \xrightarrow\, 0
\end{equation}

As $\widehat \rho(ib)=\overline{\widehat\rho(-ib)}$, we may assume without any loss of generality that $(b_n)$ is non decreasing (and positive).

\medskip
For $\varepsilon_0\in]0,\eta[$ and $n\in \N$, set 
\[
f_n:\left\{\begin{array}{rcl}B(0,\varepsilon)&\to& \C\\ z &\mapsto & |b_n|^l f(\frac z{1+|b_n|^l}+ib_n)\end{array}\right.
\]
This is a sequence of holomorphic functions on $B(0,\epsilon_0)$.

Moreover,
\[
f_n(0) =|b_n|^l (1-\widehat\rho(ib_n))\xrightarrow\, 0
\]
and the mean value inequality proves that
\[
|f_n(z)| \leqslant |f_n(0)| + |z| \frac{|b_n|^l}{1+|b_n|^l} \sup_{z\in \cal B(ib_n,\varepsilon)} |f'(z)|
\]
where the supremum is finite according to lemma~\ref{lemma:properties_Fourier_Laplace_transform}.

So, $\sup_{n\in \N} \sup_{z\in B(0,\varepsilon)}|f_n(z)|$ is finite.

This means that, $(f_n)$ is bounded in $\cal H(B(0,\epsilon_0))$ (the set of holomorphic functions on $B(0,\epsilon_0)$) and so, Montel's theorem proves that we can extract from it a subsequence that converges to an holomorphic function $f^\star$ on $B(0,\epsilon_0)$ and that the convergence is uniform on each compact subset of $B(0,\varepsilon_0)$.

But, as $\lim_{n\to +\infty} f_n(0) =0$, we have that $f^\star (0)=0$ and, as the zeros of $f^\star$ are uniformly isolated (according to lemma~\ref{lemma:uniformily_isolated_zeros}), we also get that $0$ is a simple zero of $f^\star$ and that $f^\star$ has no other $0$ in $B(0,\varepsilon_0)$.

Now, Hurwitz's theorem (see~\cite{Gam01}) proves that for any $\varepsilon<\varepsilon_0$ (we may have to decrease the value of $\varepsilon_0$), there is $N\in \N$ such that for any $n \in \N$, $n\geqslant N$, $f_n$ has exactly one zero in $B(0,\varepsilon)$.

Noting it $a_n$, and $z_n = a_n/(1+b_n^l)+ib_n$
we get that
\[
\widehat\rho(z_n) = f_n(a_n)+1 = 0 \text{ and }\Re(z_n) = \frac{\Re(a_n)}{1+|b_n|^l} \geqslant \frac{-\varepsilon}{1+|b_n|^l}
\]
This proves the lemma since for $\varepsilon_0$ small enough and $\varepsilon\in ]0,\varepsilon_0[$, we have that
\[
\frac{1+(|z_\varepsilon| + \frac{\varepsilon}{1+|b_n|})^l}{1+|z_\varepsilon|^l} \leqslant 2^l
\]
and so
\[
1+|b_n|^l \leqslant 2^l(1+|z_\varepsilon|^l)
\]
so that
\[
\frac{-\varepsilon}{1+|b_n|^l} \geqslant \frac{-2^{-l}\varepsilon }{1+|z_n|^l}
\]
\end{proof}

\section{Operators of convolution by distributions}

As we will see in section~\ref{section:speed_renewal_theorem}, the study of the speed in the renewal theorem uses the theory of operators of convolution by distributions whose Fourier transform have nice properties.

In this section, we study the following question : given a function $U$ defined on $i\R$ in $\C$, is there a distribution $T_U$ on $\R$ such that $U=\widehat{T_U}$ and what can be said about the operator of convolution by $T_U$ ?

\medskip
The first result in this direction is the following
\begin{theorem}[Paley-Wiener-Schwartz]
Let $U$ be a function on $i\R$ that can be extended to an entire function on $\C$ that satisfies that there are $M\in \N$ and $C,\tau\in\R_+^\star$ such that for any $z\in \C$,
\[
|U(z)| \leqslant C(1+|z|)^M e^{\tau|\Re(z)|}
\]
Then, $U$ is the Fourier transform of a distribution that have compact support.
\end{theorem}

The previous theorem is too restrictive for our study since it assumes that the function can be extended to an entire function whereas the functions that appear in the study of the renewal theorem are only defined on strips in $\C$ containing $i\R$. To solve this problem, we study functions that are not in $\cal D(\R)$ but who vanish exponentially fast at $\pm\infty$.

For $\eta\in\R_+^\star$ and $M\in \N$, we set
\[
\cal C^M_\eta= \left\{ f\in \cal C^M(\R) \middle| \forall k\in\lib 0,M\rib \sup_{x\in \R} e^{\eta|x|} |f^{(k)}(x)| \text{ is finite} \right\}
\]

\begin{lemma}\label{lemma:reed_simon_exponential_speed}
Let $\eta_0\in\R_+^\star$ and $U$ be an analytic function on $\C_{\eta_0}$ such that there is $M\in \N$ such that for any $\eta\in [0,\eta_0[$,
\[
\sup_{|\eta'|\leqslant \eta} \int_\R \frac{|U(\eta'+it)|}{1+|t|^M} \di t \text{ is finite}
\]
Then, for any $\eta\in ]0,\eta_0[$, $U$ is the Fourier transform of a distribution $T_U$ and the convolution by $T_U$ is continuous from $\cal C^M_\eta(\R)$ to $\cal C^0_{\eta}(\R)$.
\end{lemma}

\begin{remark}
A similar result may be proved for Fourier-Laplace transform of measures on $\Z$ to show, as an example, that any analytic function on a strip $\C_\eta$ that is also $T-$periodic is the Fourier transform of a complex measure on $\frac 1 T\Z$ whose total variation have an exponential moment. This kind of results would allow us to study the speed in the renewal theorem for lattice measures.
\end{remark}

\begin{remark}\label{remark:too_restrictive_assumption}
The integrability condition is very strong since, as an example, it is not satisfied (with $M=0$) when $U$ is the Fourier transform of a Dirac distribution (and it is clear that the convolution by Dirac distributions preserves $\cal C^0_\eta$). The problem is that in the proof, we use a theorem of Reed and Simon which proves that $U$ is the Fourier transform of a continuous function that vanishes exponentially fast when it would be more natural to say that it is the Fourier transform of a measure whose total variation has an exponential moment. Sadly, we don't know if such a theorem holds in general (i.e when $U$ is not of positive type).
\end{remark}

\begin{proof}
Let $B\in]\eta_0,+\infty[$.

Let $T_U$ be the convolution operator associated with $U$. That is to say that for any $f \in \cal S(\R)$ and any $x\in \R$,
\begin{flalign*}
T_U\star f(x) &= \frac 1 {2\pi} \int_\R e^{i\xi x} U(i\xi) \widehat f(i\xi)\di \xi =\frac 1 {2\pi} \int_\R e^{i\xi x} \frac{U(i\xi)}{(B-i\xi)^M} (B-i\xi)^M \widehat f(i\xi)\di \xi \\ 
&= \frac 1 {2\pi} \int_\R e^{i\xi x} \frac{{U(i\xi)}}{(B-i\xi)^M} \sum_{k=0}^M \binom M k B^{n-k} (-i\xi)^k \widehat f(i\xi) \di \xi \\
&=\sum_{k=0}^M \binom M k B^{n-k} \frac 1 {2\pi} \int_\R e^{i\xi x} \frac{{U(i\xi)}}{(B-i\xi)^M}  \widehat{f^{(k)}}(i\xi) \di \xi \\
\end{flalign*}

But, the function $U_0:z\mapsto \frac{U(z)}{(B-z)^M}$ is analytic on $\C_{\eta_0}$ and satisfies that for any $\eta\in]0,\eta_0[$,
\[
\sup_{|\eta'| \leqslant \eta} \int_\R |U_0(\eta'+it)| \di t\text{ is finite}
\]
So, according to Reed and Simon~\cite[Theorem IX.14]{RS75}, there is a function $g\in \cal C^0(\R)$ such that for any $\eta\in ]0,\eta_0[$,
\[
\sup_{x\in \R} e^{\eta|x|} |g(x)| \text{ is finite}
\]
And, $\widehat g(z) = U_0(z)$.

Thus, for any $k\in [0,M]$,
\[
\frac 1 {2\pi}\int_\R e^{i\xi x} \widehat g(i\xi) \widehat{f^{(k)}}(\xi)\di \xi = \int_\R g(t) f^{(k)}(x-t) \di t
\]
and finally, we get that for any $f\in \cal S(\R)$,
\[
T_U\star f(x) =  \sum_{k=0}^M \binom M k B^{M-k} \int_\R g(t) f^{(k)}(x-t) \di t
\]
Now, since $\cal S(\R)$ is dense in $\cal C^M_\eta$, this formula may be extended to $\cal C^M_\eta(\R)$ and this defines the wanted operator on $\cal C^M_\eta(\R)$.

Finally, for any $f\in \cal C^M_\eta(\R)$, $T_U\star f$ is continuous since it is the sum of convolutions of continuous functions, and $T_U\star f\in \cal C^0_{\eta}(\R)$ since, for any $h\in \cal C^0_\eta(\R)$,
\[
e^{\eta|x|} \left| \int_\R g(t) h(x-t)\di t \right| \leqslant e^{\eta |x|} \int_\R |g(t)| e^{-\eta|x-t|} \di t \leqslant \int_\R |g(t)|e^{\eta|t|}\di t \text{ which is finite}
\]
where we used that for any $x,t\in \R$ and any $\eta \in \R_+^\star$,
\[
e^{\eta(|x| - |x-t|)} \leqslant e^{\eta ||x| - |x-t||} \leqslant e^{\eta|t|} 
\]
\end{proof}

The following definition and proposition are the translation of the compactly supported distributions and Paley-Wiener-Schwartz's theorem for tempered distributions.

\begin{definition}[Tempered distributions of fast decrease]
Let $T\in \cal S'(\R)$.

We say that $T$ is of fast decrease if for any $f\in \cal S(\R)$ we have that $T\star f\in \cal S(\R)$.
\end{definition}

\begin{proposition}\label{lemma:multiplier_Schwartz_space}
Let $U\in \cal C^\infty(\R)$ be such that for any $n\in \N$, $t\mapsto |U^{(n)}(t)|$ is at most polynomial.

Then, $U$ is the Fourier transform of a tempered distribution of fast decrease.
\end{proposition}

\subsection{Weighted Sobolev spaces}
Until now, we studied the problem in Frechet spaces. To work with Banach spaces, we will use weighted Sobolev spaces.

\begin{definition}[Weight]
We call weight any continuous function $\omega:\R \to [1,+\infty[$.
\end{definition}

From now on, we fix a weight $\omega$.

\medskip
For $m\in \N$ and $p\in[1,+\infty[$, we set
\[
\cal H^{m,p}_\omega(\R) = \left\{ f\in \rm L^2(\R)\middle| \forall k\in\lib 0,\dots,m\} \;\;\int_\R \omega(x) |f^{(k)}(x)|^p \di x \text{ is finite} \right\}
\]
and
\[
\cal H^{m,\infty}_\omega(\R) = \left\{ f\in \rm L^2(\R)\middle| \forall k\in\{0,\dots,m\}\;\;\omega f^{(k)} \in \rm{L}^\infty(\R)\right\}
\]
If $m=0$ and $p\in[1,\infty]$, we note
\[
\mathrm{L}^p_\omega(\R)=\cal H^{m,p}_\omega(\R) 
\]

These spaces are Banach spaces are Banach spaces when they are endowed with the natural norms
\[
\|f\|_{\cal H^{m,p}_\omega(\R)} = \max_{k\in\{0,\dots, m\}} \left( \int_\R \omega(x) |f^{(k)}(x)|^p \di x \right)^{1/p}
\]
and
\[
\|f\|_{\cal H^{m,\infty}_\omega(\R)} = \max_{k\in\{0,\dots, m\}} \|\omega f^{(k)} \|_\infty
\]

Moreover, for any $m,l\in \N$ and any $p\in [1,+\infty]$,
\[
\cal C^\infty_c(\R) \subset \cal H^{m+l,p}_\omega(\R)  \subset \cal H^{m,p}_\omega(\R) 
\]
and the inclusion are continuous.

Thus, if there are no Fourier multipliers $P:\cal C^\infty_c(\R) \to \mathrm{L}^p_\omega(\R)$, then there are no multipliers $P:\cal C^\infty_c(\R) \to \cal H^{m,p}_\omega(\R) $.

\medskip
To begin our study of operators of convolutions on weighted Sobolev spaces, we give a result of L\"ofstr\"om stated in~\cite{Lof83}.
\begin{theorem*}
For $A\in \R_+^\star$, $\alpha\in]1,+\infty|$ and $x\in \R$, let $\omega(x) = e^{A|x|^\alpha}$, then there is no non-trivial continuous Fourier multiplier on $\rm{L}^p_\omega(\R)$. More precisely, if $P$ is such an operator, then $P=a\delta_0$ for some $a\in\R$.
\end{theorem*}

Thus, we will no longer study weights such as in the previous theorem and make the following
\begin{definition}
A weight $\omega$ is called sub-multiplicative or a Beurling weight if for any $x,y\in \R$,
\[
\omega(x+y)\leqslant \omega(x) \omega(y)
\]
$\omega$ is called sub-exponential if
\[
\sup_{x\in \R} \frac{\ln(\omega(x))}{|x|} \text{ is finite}
\]

Moreover, it is called strictly sub-exponential if for any $\varepsilon\in\R_+^\star$,
\[
\lim_{x\to \pm\infty} e^{-\varepsilon |x|} \omega(x) = 0
\]
\end{definition}

\begin{remark}
If for some $C\in \R_+$ we have that for any $x,y\in \R$, $\omega(x+y) \leqslant C\omega(x) \omega(y)$ then define $\tilde \omega(x) = C\omega(x)$ and note that $\tilde\omega$ is sub-multiplicative and that the spaces we are going to define for $\omega$ are isomorphic to the ones we get taking $\tilde\omega$ instead.
\end{remark}

\begin{example}~
\begin{itemize}
\item For any $A\in \R_+^\star$ and $\alpha \in]1,+\infty[$, $\omega(x) = e^{A|x|^\alpha}$ is not a sub-multiplicative.
\item For any $A\in \R_+$, $\omega(x) = e^{A|x|}$ is sub-multiplicative.
\item For any $\alpha\in[0,1[$ and $A\in \R_+$, $\omega(x) = e^{A|x|^\alpha}$ is strictly sub-multiplicative.
\item For any $k\in \N$, $\omega(x) = (1+|x|)^k$ is strictly sub-multiplicative and we will call it polynomial.
\end{itemize}
\end{example}

\begin{remark}
Young's inequality proves that if $\omega$ is sub-multiplicative then for any $p\in [1,+\infty]$ $\mathrm{L}^1_\omega \cap\mathrm{L}^p_\omega$ is a Banach algebra.
\end{remark}

Let $P:\cal H^{m+l,p}_\omega(\R)  \to \cal H^{m,p}_\omega(\R) $ be a continuous operator of convolution by some distribution whose Fourier-transform is denoted by $U$. Saying that $P$ is continuous means that for any $f\in \cal H^{m+l,p}_\omega(\R) $,
\[
\int_\R \omega(x) |Pf^{(k)}(x)|^p \di x = \int_\R \omega(x) \left| \cal F(U\widehat f)^{(k)}(x)\right|^p \di x \text{ is finite}
\]
and this links $\omega$ and $U$.

In this subsection, we intend to study this. First of all, we will assume that $\omega$ is polynomial, then we will study what happens if $U$ can be extended to a meromorphic function on a neighbourhood of the imaginary axis with poles in a particular domain and show that we cannot have a speed that is faster than polynomial. Then, we will study the case when $U$ can be extended to an analytic function on a strip and see that, as in lemma~\ref{lemma:reed_simon_exponential_speed}, we can take an exponential weight.

\subsection{Polynomial speed}
First of all, the following lemma makes precise the idea that the Fourier transform inverts the regularity of a function and it's vanishing at infinity.

\begin{lemma}[see {\cite[Theorem 1.1]{Bou04}}]\label{lemma:fourier_polynomial_sobolev}
For any $n\in \N$, and any $x\in \R$, note $p_n(x) = 1+|x|^n$.

Then, for any $k,m\in \N$, the Fourier transform is a continuous isomorphism between $\cal H^{m,2}_{p_k}(\R)$ and $\cal H^{k,2}_{p_m}(\R)$.
\end{lemma}

The previous lemma has the direct following

\begin{corollary}\label{corollary:polynomial_sobolev}
Let $U\in \cal C^n(i\R)$ a function on $i\R$ such that for any $p\in [0,n]$,
\[
U^{(p)}(it) \in \cal O(1+|t|^{l})
\]
Then, for any $m,k\in \N$, $k\leqslant n$, the Fourier multiplier defined by $U$ is continuous from $\cal H^{m+2l,2}_{p_{k}}(\R) $ onto $\cal H^{m,2}_{p_k}(\R)$.
\end{corollary}

\begin{proof}
Having lemma~\ref{lemma:fourier_polynomial_sobolev}, we only need to prove that the multiplication by $U$ is continuous from $\cal H^{k,2}_{p_{m+2l}}(\R)$ onto $\cal H^{k,2}_{p_m}(\R)$.

Note
\[
C_U=\sup_{p\in[0,n]} \sup_{x\in \R} \frac{\left|U^{(p)}(ix)\right|}{(1+|x|)^l} 
\]

For any $r\in [0,k]$ and any $f\in \cal H^{k,2}_{p_{m+2l}}(i\R)$,
\begin{flalign*}
\int_\R (1+|x|)^{m} \left| \left(fU\right)^{(r)} (ix) \right|^2 \di x &\leqslant \sum_{p=0}^r \binom r p \int_\R (1+|x|)^{m} \left| f^{(p)}(ix)U^{(r-p)} (ix) \right|^2 \di x \\
& \leqslant \sum_{p=0}^k \binom k p C_U \int_\R (1+|x|)^{m+2l} \left|f^{(p)}(ix) \right|^2 \di x
\end{flalign*}
Thus, $fU \in \cal H^{k,2}_{p_m}(\R)$ and the multiplication by $U$ is continuous.
\end{proof}

\begin{lemma}
The Fourier transform maps $\cal H_{p_m}^{k,1}$ to $\cal H_{p_k}^{m,\infty}$ and $\cal H_{p_m}^{k,\infty}$ to $\cal H_{p_k}^{m-2,\infty}$.
\end{lemma}

\begin{proof}
\begin{flalign*}
\left|\xi^l \widehat f^{(n)}(\xi) \right|  \leqslant \sum_{r=0}^{\min(n,l)} \binom l r \frac{n!}{(n-r)!} \int_\R |x|^{n-r} \left| f^{(l-r)} (x) \right| \di x 
\end{flalign*}

But, if $f\in \cal H_{p_m}^{k,1}$, then for any $l\leqslant k$, any $n\leqslant m$ and any $r\leqslant \min(l,n)$, we have that
\[
\int_\R (1+|x|^{n-r}) \left| f^{(l-r)} (x) \right| \di x  \leqslant \|f\|_{ \cal H_{p_m}^{k,1}}
\]
and this proves the first part of the lemma.

\medskip
To prove the second one, note that if $f\in \cal H_{p_m}^{k,\infty}$ we have that for any $l\leqslant k$,
\[
\left|f^{(l)}(x) \right| \leqslant \|f\|_{ \cal H_{p_m}^{k,\infty}} \frac 1 {1+|x|^m}
\]
So, for any $l\leqslant k$ and any $n\in \N$,
\[
\left|\xi^l \widehat f^{(n)}(\xi) \right|  \leqslant \sum_{r=0}^{\min(n,l)} \binom l r \frac{n!}{(n-r)!} \int_\R |x|^{n-r} \frac {\|f\|_{ \cal H_{p_m}^{k,\infty}}} {1+|x|^m} \di x 
\]
And so, $\widehat f \in \cal H^{m-2,\infty}_{p_k}$.
\end{proof}

As in corollary~\ref{corollary:polynomial_sobolev}, the previous lemma has next
\begin{corollary}\label{corollary:polynomial_conergence_norme_infinity}
Let $l,n\in \N$ and $U\in \cal C^n(i\R)$ be a function on $i\R$ such that for any $p\in [0,n]$,
\[
U^{(p)}(it) \in \cal O(1+|t|^{l})
\]
Then, for any $m,k\in \N$, $k\leqslant n$, the convolution by $T_U$ is continuous from $\cal H^{m+2+l,\infty}_{p_{k+2}}(\R) $ onto $\cal H^{m,\infty}_{p_k}(\R)$.
\end{corollary}

\subsection{Exponential speed}

For $a\in\R_+^\star$, we note
\[
\mathcal{H}^{2,m}_{a}(\R)=\left\{ f\in \mathrm{L}^2(\R) \middle| \forall k\in[0,m]  \;e^{a|x|}f^{(k)}(x) \in \mathrm{L}^2(\R)\right\}
\]

First, we extend a result of Harper in~\cite{Har09} that gives the range of the Fourier transform of $\mathcal{H}^{2,0}_a(\R)$ in next 
\begin{lemma}\label{lemma:paley_wiener_exponential}For any $a\in \R_+^\star$ and any $m\in \N$, the Fourier transform is an isomorphism between $\cal H^{2,m}_a(\R)$ and
\[
\cal F\cal H^{2,m}_a(\R)=\left\{ f\in \cal H(\C_a) \middle| \int_{[-a,a]\times \R}(1+|x+iy|^m)^2|f(x+iy)|^2 \di x\di y \text{ is finite} \right\}
\]
\end{lemma}

\begin{proof}
Let $f\in \left\{ f\in \cal H(\C_a) \middle| \int_{[-a,a]\times \R}(1+|x+iy|^m)^2|f(x+iy)|^2 \di x\di y \text{ is finite}\right\}$.

For any $k\in [0,m]$, note $f_k(z) = (-z)^k f(z)$, we have that $f_k$ is analytic on $\C_a$ and that $\int_{[-a,a]\times \R} |f_k(x+iy)|^2 \di x\di y$ is finite. So, according to~\cite[Theorem 2.1]{Har09}, there is $g_k\in \mathrm{L}^2(\R)$ such that $(x\mapsto e^{a|x|} g_k(x) )\in \mathrm{L}^2(\R)$ and $f_k = \widehat g_k$. 

Thus, if we manage to prove that for any $k\in [1,m]$, $f_k= f_0^{(k)}$ the first part of the lemma will be proved since this will imply that $g_0 \in \mathcal{H}^{2,m}_{a}(\R)$ and we also have that $\widehat g_0=f$.

But, for any $\varphi \in \cal S(\R)$,
\begin{flalign*}
\int_\R f_0^{(k)} (x) \varphi(x) \di x &= (-1)^k \int_\R f_0(x) \varphi^{(k)} (x) \di x = (-1)^k \int_\R \widehat f_0(i\xi) \overline{ \varphi^{(k)}(i\xi) }\di \xi \\
&=(-1)^k \int_\R \widehat f_0(i\xi) \overline{(-i\xi)^k \widehat \varphi(i\xi)} \di \xi = \int_\R f_k(i\xi) \overline{\widehat \varphi(i\xi)} \di \xi = \int_\R f_k(x)\varphi(x) \di x 
\end{flalign*}
and this proves that for any $k\in[0,m]$, $f_k = f_0^{(k)}$.

\medskip
Take now $f\in \mathcal{H}^{2,m}_{a}(\R)$. Then, \cite[Theorem 2.1]{Har09} proves that for any $k\in [0,m]$, $\widehat {f^{(k)}}\in \cal H(\C_a)$ and $\int_{[-a,a] \times \R} |\widehat{f^{(k)}}(x+iy)|^2 \di x \di y$ is finite.

But, for any $z\in \C_a$, $\widehat{f^{(k)}}(z) = z^k \widehat f(z)$ and so, we get that
\[
\int_{[-a,a] \times \R}(1+|x+iy|^{2m}) |\widehat{f}(x+iy)|^2 \di x \di y
\]
We conclude since there is $C\in R$ such that for any $z\in \C_a$,
\[
1+|z|^{2m} \leqslant C(1+|z|^m)^2
\]
\end{proof}

\begin{corollary}\label{corollary:exponential_sobolev}
Let $a\in \R_+^\star$ and $U$ be an analytic function on $\C_a$ such that there is $l\in \N$ such that $U(z) \in \cal O(1+|z|^l)$.

Then, for any $m\in\R$, $U$ is the symbol of a Fourier multiplier $P: \cal H^{2,m+l}_a(\R) \to \cal H^{2,m}_a(\R)$.
\end{corollary}

\begin{proof}
The condition on $U$ implies that the multiplication by $U$ continuously maps $\cal F\cal H^{2,m+l}_a(\R)$ onto $\cal F\cal H^{2,m}_a(\R)$. So the corollary is a direct consequence of lemma~\ref{lemma:paley_wiener_exponential}.
\end{proof}

\subsection{Between polynomial and exponential speed}
The aim of this subsection is to study what happens in the intermediate case when the function $U$ has a meromorphic extension to a strip containing $i\R$ but there is a sequence of poles of $U$ whose real part converges to $0$.

We will not study this problem in such a generality but we will assume that we have some control on $U$ and this will actually be enough to study the renewal theorem. The following definition is technical but the reader may think that it is exactly what the Fourier transform of $G-T_\lambda$ will satisfy (see the introduction for a definition of $G$ and $T_\lambda$).

\begin{definition}\label{definition:L_eta}Let $\eta \in \R_+^\star$.
We note $\cal L_\eta$ the set of meromorphic functions $U$ on $\C_\eta$ that satisfies the following assumptions
\begin{enumerate}
\item\label{item:simple_poles} The poles of $U$ are simple
\item\label{item:separated_poles} The set $\cal A$ of poles of $U$ is infinite  and separated.
\item\label{item:shape_poles_U} There are $l_1,l_2\in \R_+^\star$ such that
\[
\sup_{a\in \cal A}{(1+|a|)^{l_1}} \Re(a)  <0 \text{ and }\sup_{a\in \cal A}{(1+|a|)^{l_2}} \Re(a)=0
\]
\item\label{item:residues} There are $C_0,l_0 \in \R_+$ such that for any $a\in \cal A$
\[
\frac 1 {C_0 (1+|a|)^{l_0}} \leqslant |Res_a(U)| \leqslant C_0(1+|a|)^{l_0}
\]
\item\label{item:control_outside} There are $\varepsilon\in \R_+^\star$ and $l\in \R$ such that
\[
\sup_{\substack{z\in \C_\eta\\d(z,A)\geqslant \varepsilon}} \frac{|U(z)|}{(1+|z|)^l} \text{ is finite}
\]
\end{enumerate}
\end{definition}

\begin{remark}
Note that under these assumptions we have that for any $a\in \cal A$, $\Re(a)<0$.
\end{remark}

In next lemma, we give an equivalent definition of $\cal L_\eta$ that is easier to deal with.

\begin{lemma}\label{lemma:other_definition_L_eta}
The set $\cal L_\eta$ is the set of functions that writes
\[
U(z) = U_0(z) + \sum_{a\in \cal A} \left(\frac{z}{a}\right)^{m} \frac {\lambda_a}{z-a}
\]
where $\cal A$ is a separated set such that for some $l_1,l_2\in \R_+^\star$,
\[
\sup_{a\in \cal A}{(1+|a|)^{l_1}} \Re(a)  <0 \text{ and }\sup_{a\in \cal A}{(1+|a|)^{l_2}} \Re(a)=0
\]
The function $U_0$ is holomorphic on $\C_\eta$ such that for some $l\in \R$, $U_0(z) \in \cal O((1+|z|)^l)$ and $(\lambda_a)\in \C^{\cal A}$ is such that there are $C_0,l_0\in  \R_+$ such that for any $a\in \cal A$,
\[
\frac 1 {C_0 (1+|a|)^{l_0}} \leqslant |\lambda_a| \leqslant C_0(1+|a|)^{l_0}
\]
and $m\geqslant l_0+1$.
\end{lemma}

\begin{proof}
Let $U\in \cal L_\eta$ and note $\cal A$ it set of poles.

Then,
\[
U_1(z) = \sum_{a\in \cal A} \left( \frac z a \right)^{l_0+1} \frac {Res_a(U)}{z-a}
\]
defines a meromorphic function and $U-U_1$ can be extended to an holomorphic function on $\C_\eta$. So, we only need to prove that for some $l_1\in \R_+$,
\[
\sup_{z\in \C_\eta} \frac{|U(z)-U_1(z)|}{1+|z|^{l_1}} \text{ is finite}
\]
But, this is just the maximum principle since
\[
\sup_{\substack{z\in \C_\eta\\ d(z,\cal A) \geqslant \varepsilon}} \frac{|U(z)-U_1(z)|}{1+|z|^{l_1}}  \leqslant \sup_{\substack{z\in \C_\eta\\ d(z,\cal A) \geqslant \varepsilon}} \frac{|U(z)|}{1+|z|^{l_1}} + \frac 1 \varepsilon \sup_{z\in \C_\eta} \frac{|z|^{l_0+1}}{1+|z|^{l_1}} \sum_{a\in \cal A} \frac{|Res_a(U)|}{|a|^{l_0+1}}
\]
so, if we take $l_1$ large enough, we get the expected result. 
\end{proof}

As the assumptions on $U\in \cal L_\eta$ are between the ones of lemma~\ref{lemma:reed_simon_exponential_speed} and \ref{corollary:polynomial_conergence_norme_infinity}, we wonder if for any $f\in \cal C^{\infty}_c(\R)$, $T_U\star f$ vanishes very fast at $\pm\infty$ (maybe at a speed between polynomial and exponential since we already know, with lemma~\ref{lemma:multiplier_Schwartz_space} that it vanishes faster than any polynomial). This is actually false in general as we will prove in next
\begin{proposition}\label{proposition:between_polynomial_exponential}
Let $U\in \cal L_\eta$ and $\omega\in \Omega$ (see definition~\ref{definition:Omega_0}).

Assume that the convolution by $T_U$ is continuous from $\cal C^\infty_c(\R)$ to $\mathrm{L}_\omega^2(\R)$ then $\omega$ doesn't grow faster than any polynomial.
\end{proposition}

\begin{remark}
The condition on $\omega \in \Omega$ only deals with the one-sided Fourier-Laplace transform of $\omega$. This is because we assumed that the poles of $U$ are in $\C_-$.
\end{remark}

\begin{remark}
The functions on $\R$ that we are studying are implicitly supposed to take complex values thus we don't need to ask symmetry conditions on $\cal A$ and on $(\lambda_a)$ (this means that we may have some $a\in \cal A$ such that $\overline a \not\in \cal A$ and some $a\in \cal A$ such that $\overline{a} \in \cal A$ but $\lambda_{\overline a} \not= \overline{ \lambda_a}$).
\end{remark}

Next lemma is the first step of the proof of proposition~\ref{proposition:between_polynomial_exponential}. It shows that if a sum of complex exponential functions belong to $\mathrm{L}^2_\omega(\R_+)$ then it is not because of cancellations since the sum of the squared modulus also belong to $\mathrm{L}^2_\omega(\R_+)$. This is possible because we defined $\Omega$ to be the set of functions whose Fourier-Laplace transform is bounded close to the imaginary axis except at $0$.
\begin{lemma}\label{lemma:control_theta}
Let $\eta \in \R_+^\star$ and $\cal A\in ]-\eta,0[ \oplus i\R$ be such that
\[
\inf_{\substack{a_1,a_2\in\cal A\\ a_1\not=a_2}} |a_1-a_2|=\delta>\eta
\]
Let $(u_a) \in \ell^1(\cal A)$.

Then, for any $\omega\in \Omega$ (see definition~\ref{definition:Omega_0} and equation~\ref{equation:theta_omega} for the definition of $\Omega$ and $\Theta_\omega$)
\[
\sum_{a\in \cal A} |u_a|^2 \int_{\R_+} e^{2\Re(a)x} \omega(x) \di x \leqslant \int_{\R_+} \omega(x)\left|\sum_{a\in \cal A} u_a e^{ax} \right|^2 \di x + \Theta_\omega(\delta-\eta) \left(\sum_{a\in \cal A} |u_{a}  | \right)^2
\]
(all the integrals may be infinite)
\end{lemma}

\begin{proof}
For any $x\in \R_+$, we have that
\[
\left|\sum_{a\in \cal A} u_a e^{ax} \right|^2 = \sum_{a\in \cal A} |u_a|^2 e^{2\Re(a)x} + \sum_{a_1,a_2\in \cal A,\;a_1\not=a_2} u_{a_1} \overline{u_{a_2}} e^{(a_1+\overline{a_2})x}
\]
where all the involved series are absolutely convergent since for any $a\in \cal A$, $\Re(a)<0$ and $(u_a)\in \ell^1(\cal A)$.

Moreover,
\[
\left|\sum_{a_1,a_2\in \cal A,\;a_1\not=a_2} u_{a_1} \overline{u_{a_2}} e^{(a_1+\overline{a_2})x}\right| \leqslant \sum_{a_1,a_2\in \cal A,\;a_1\not=a_2} |u_{a_1} \overline{u_{a_2}}|
\]
so for any $s\in \R_+^\star$, the function
\[
x\mapsto \omega(x)e^{-sx} u_{a_1} \overline{u_{a_2}} e^{(a_1+\overline{a_2})x} 
\]
is absolutely integrable since $\omega$ is strictly sub-exponential and this means, using Fubini's theorem, that
\[
\int_{\R_+} \omega(x) e^{-sx}\sum_{a_1,a_2\in \cal A,\;a_1\not=a_2} u_{a_1} \overline{u_{a_2}} e^{(a_1+\overline{a_2})x} \di x= \sum_{\underset{a_1\not=a_2}{a_1,a_2\in \cal A}} u_{a_1} \overline{u_{a_2}} \int_{\R_+} \omega(x)e^{-sx} e^{(a_1+\overline{a_2})x} \di x
\]

Therefore, for any $s\in \R_+^\star$,
\begin{flalign*}
\sum_{a\in \cal A} |u_a|^2 \int_{\R_+} \omega(x) e^{-sx} e^{2\Re(a)x} \di x &\leqslant \int_{\R_+} \omega(x)e^{-sx} \left|\sum_{a\in \cal A} u_a e^{ax} \right|^2 \di x \\& \retrait + \sum_{\underset{a_1\not=a_2}{a_1,a_2\in \cal A}} |u_{a_1} \overline{u_{a_2}} | \left|\int_{\R_+} \omega(x) e^{-sx} e^{(a_1+\overline{a_2})x} \di x \right|
\end{flalign*}
But, $\cal A$ is $\delta-$separated so if $a_1,a_2\in \cal A$ are such that $a_1\not=a_2$, then $|a_1-a_2|\geqslant \delta$, thus 
\[
|a_1+\overline{a_2}-s|\geqslant|\Im(a_1-a_2)|\geqslant |a_1-a_2| - |\Re(a_1-a_2)| \geqslant \delta-\eta
\]
and so, by definition of $\Theta_\omega$,
\[
\left|\int_{\R_+} \omega(x) e^{-sx}e^{(a_1+\overline{a_2})x} \di x \right| \leqslant \Theta_\omega(\delta-\eta)
\]
thus, for any $s\in \R_+^\star$,
\begin{flalign*}
\sum_{a\in \cal A} |u_a|^2 \int_{\R_+} \omega(x)e^{-sx} e^{2\Re(a)x} \di x & \leqslant \int_{\R_+}  \omega(x) e^{-sx}\left|\sum_{a\in \cal A} u_a e^{ax} \right|^2 \di x \\
& \retrait + \Theta_\omega(\delta-\eta)\sum_{\underset{a_1\not=a_2}{a_1,a_2\in \cal A}} |u_{a_1} \overline{u_{a_2}} |
\end{flalign*}
and we conclude using the monotone convergence theorem.
\end{proof}

The following lemma is the main step in the proof of proposition~\ref{proposition:between_polynomial_exponential}.
\begin{lemma}\label{lemma:between_speeds}
Let $\omega\in \Omega$ and $U\in \cal L_\eta$.

Assume that the convolution by $T_U$ is continuous from $\cal C^\infty_c(\R)$ to $\mathrm{L}_{\omega}^2(\R)$.

Then, there is $p\in \N$ such that
\[
\int_{\R_+} \omega(x) e^{2\Re(a)x} \di x  \in \cal O({1+|a|^p})
\]
\end{lemma}

\begin{remark}
Before the proof, remark that we obtain a condition on $\omega$ that only depend on the values of $\omega$ on $\R_+$. This is because we assumed that all the poles of $U$ are in $\C_-$.
\end{remark}

\begin{proof}
Let $U\in \cal L_\eta$ and write the decomposition given by lemma~\ref{lemma:other_definition_L_eta}
\[
U(z) = U_0(z) + \sum_{a\in \cal A} \left(\frac{z}{a}\right)^{m} \frac {\lambda_a}{z-a}
\]
The operator $T_U$ of convolution by the inverse Fourier-Laplace transform of $U$ writes, $f\in \cal C^\infty_c(\R)$ and $x\in \R$,
\begin{flalign*}
T_U\star f(x) &= \frac 1 {2\pi}\int_\R \widehat f(\xi) e^{i\xi x} U_0(i\xi) \di \xi +\frac 1 {2\pi} \int_\R \widehat f(\xi) e^{i\xi x} \sum_{a\in \cal A}\left(\frac{i\xi}a\right)^m \frac{ \lambda_a}{i\xi-a} \di \xi \\
&=\frac 1 {2\pi}\int_\R \widehat f(\xi) e^{i\xi x} U_0(i\xi) \di \xi + \sum_{a\in\cal A} \frac {\lambda_a}{a^m} \int_\R f^{(m)}(x-u) \un_{\R_+}(u) e^{au} \di u \\
&= \frac 1 {2\pi} \int_\R \widehat f(\xi) e^{i\xi x} U_0(i\xi) \di \xi + \sum_{a \in \cal A}\frac {\lambda_a}{a^m} e^{ax} \int_{-\infty}^x f^{(m)}(u)e^{-au} \di u
\end{flalign*}
Note
\[
g(x)= \frac 1 {2\pi} \int_\R \widehat f(\xi) e^{i\xi x} U_0(i\xi) \di \xi
\]
then, according to lemma~\ref{lemma:reed_simon_exponential_speed}, $g\in \cal C^0_\eta$ since $U_0(z) \in \cal O(1+|z|^l)$ so, $g\in\mathrm{L}^2_\omega(\R)$ since $\omega$ is sub-exponential and the application that maps $f$ to $g$ is continuous.

Thus, if $T_U:\cal C^\infty_c(\R)\to \mathrm{L}^2_\omega(\R)$ is continuous, then so is the operator that maps $f$ to $Pf-g$ and this means that for any $M\in \R_+^\star$ there are $C_M,N_M$ such that for any $f\in \cal C^\infty([-M,M])$,
\[
\int_\R \omega(x) \left| \sum_{a \in \cal A}\frac {\lambda_a}{a^m} e^{ax} \int_{-\infty}^x f^{(m)}(u)e^{-au} \di u \right|^2 \di x \leqslant C_M^2\max_{k\in[0,N_M]} \|f^{(k)}\|_\infty^2
\]
Moreover, 
\[
\int_{-\infty}^x f^{(m)}(u)e^{-au}\di u = \widehat{f^{(m)}}(a) - \int_x^{+\infty} f^{(m)}(u) e^{-au} \di u
\]
and so, for $x\geqslant M$, we obviously have that $\int_{x}^{+\infty} f^{(m)}(u) e^{-au}\di u=0$.

Thus, (this is where we pass from $\R$ to $\R_+$)
\[
\int_{\R_+}\omega(x)\left|\sum_{a \in \cal A}\frac{\lambda_a}{a^m} e^{ax} \int_{x}^{+\infty} f'^{(m)}(u)e^{-au} \di u \right|^2 \di x \leqslant \sum_{a\in \cal A} \frac{|\lambda_a|^2}{|a|^{2m}} M^2\|f''\|_\infty^2\int_{0}^M \omega(x) \di x
\]
And this finally proves, for some other constant that we re-note $C_M$ that only depend on $\cal A$, $(\lambda_a)$, $M$ and $\omega$, we have that for any $f \in \cal C^\infty_c([-M,M])$,
\[
\int_{\R_+} \omega(x) \left| \sum_{a \in \cal A}\lambda_a e^{ax}\widehat f(a) \right|^2 \di x \leqslant C_M^2 \max_{k\in [0,N_M]} \|f^{(k)}\|_\infty^2\text{ (we used that }\widehat{f''}(a) = a^2 \widehat f(a)\text{)}
\]
Using lemma~\ref{lemma:control_theta} (whose assumptions hold since $\widehat f(a)$ decreases faster than any polynomial for $f\in \cal C^\infty([-M,M])$ according to lemma~\ref{proposition:strip_control_Laplace_dirac}) we get that for any $M$ in $\R_+^\star$, there are $C_M$ and $N_M$ such that for any $f\in \cal C^\infty([-M,M])$,
\[
\sum_{a\in \cal A} |\lambda_a|^2 |\widehat f(a)|^2 \int_{\R_+} \omega(x) e^{2\Re(a)x} \di x\leqslant C_M^2 \max_{k\in [0,N_M]} \|f^{(k)}\|_\infty^2
\]
and, this precisely implies that $(w_a \widehat \delta_a)$ is bounded in $\cal D'(\R)$ where
\[
w_a = |\lambda_a| \left(\int_{\R_+} \omega(x) e^{2\Re(a)x} \di x \right)^{1/2}
\]
And so, proposition~\ref{proposition:strip_control_Laplace_dirac} proves that there is $p\in \R_+$ such that
\[
|\lambda_a|^2  \int_{\R_+} \omega(x) e^{2\Re(a)x} \di x  \in \cal O({1+|a|^p})
\]
and, as we assumed in the definition of $(\lambda_a)$ that
\[
\sup_{a\in \cal A}{(1+|a|)^{l_1}} \Re(a)  <0 
\]
we have that for some $C_\lambda\in \R_+^\star$, and any $a\in \cal A$,
\[
|\lambda_a| \geqslant \frac{C_\lambda}{(1+|a|)^{l_1}}
\]
and so,
\[
\int_{\R_+} \omega(x) e^{2\Re(a)x} \di x \in \cal O(1+|a|^{p+2l_1})
\]
which is what we intended to prove.

\end{proof}

\begin{proof}[End of the proof of proposition~\ref{proposition:between_polynomial_exponential}]
We proved in lemma~\ref{lemma:between_speeds} that under the assumptions of the proposition, there is $p\in \R_+$ such that
\[
  \int_{\R_+} \omega(x) e^{2\Re(a)x} \di x  \in \cal O({1+|a|^p})
\]
But, remember $\omega$ is non negative and that, by assumption, there are $C_2,l_2\in \R_+^\star$ such that for infinitely many $a\in \cal A$,
\[
\frac{-C}{1+|a|^{l_2}} \leqslant \Re(a)<0
\]
So, we obtain that
\[
\int_{\R_+} \omega(x) e^{2\Re(a)x} \di x  \in \cal O(|\Re(a)|^{-p/l_2})
\]
So, as $\sup_{a\in \cal A}\Re(a)=0$, lemma~\ref{lemma:essential_singularity_Laplace_transform} proves that $\omega$ doesn't grow faster than any polynomial.
\end{proof}

\section{The speed in the renewal theorem}\label{section:speed_renewal_theorem}

In this section, we finally use what we did on weakly-diophantine measures and Fourier multipliers to study the speed in Kesten's renewal theorem.

\subsection{The renewal theorem}

Let $\rho$ be a probability measure on $\R$ which have an exponential moment and a positive drift $\lambda=\int y\di\rho(y)>0$.

We call ``Green kernel'' the measure
\begin{equation}\label{equation:G}
G  = \sum_{n=0}^{+\infty} ( \widetilde{\rho})^{\star n} \text{ where for any borelian subset }A\text{ of }\R,\; \widetilde\rho(A)= \int_\R \un_A(-y) \di\rho(y)
\end{equation}
The large deviation inequality proves that for any $x\in \R$, $G\star \delta_x$ is finite on borelian bounded subsets of $\R$.

\medskip
\begin{remark}
One has to be careful to the inversion of $\rho$ in the definition of $G$ which is made to have the following equation for any $f\in \cal C^0_c(\R)$ and any $x\in \R$,
\[
G\star f(x) = \sum_{n=0}^{+\infty} P^nf(x)
\]
Where $P$ is the Markov operator associated to $\rho$ defined by $Pf(x) = \int_\R f(x+y) \di\rho(y)$.
\end{remark}

To understand the measure $G$, we first have the renewal theorem (see in~\cite{Fel71} Chapter XI, the ninth-section) which tells us that, if $\rho$ is non lattice (see definition~\ref{definition:lattice_measure}), then for any $f\in\cal C^0_c(\R)$,
\[
\lim_{x\to -\infty } <G\star \delta_x,f > = \frac {1} \lambda \int_\R f(u)\di u \text{ and }\lim_{x\to +\infty} <G\star \delta_x,f> = 0
\]

We define a measure $T_\lambda$ on $\R$ by
\begin{equation}\label{equation:T_lambda}
\forall f\in \cal C^0_c(\R) \;\;<T_\lambda ,f> = \frac{1} \lambda \int_{\R_-} f(u) \di u 
\end{equation}

Thus,
\[
T_\lambda \star f(x) = <T_\lambda \star \delta_x ,f > =\frac{1} \lambda \int_{\R_-} f(x-u) \di u= \frac{1} \lambda \int_{x}^{+\infty} f(u)\di u
\]
And, with these notations, the renewal theorem becomes :
\begin{equation}\label{equation:renewal_theorem}
\text{ in }\cal C^0_c(\R)^\star ,\;(G-T_\lambda)\star \delta_x  \displaystyle{\smash{\,\mathop{\rightharpoonup}\limits_{x\to\pm\infty}^{\star}\,}} 0
\end{equation}
Moreover, if $f$ converges to $0$ polynomially fast (resp. exponentially fast) at $\pm\infty$, then $T_\lambda \star f$ converges to it's limits at almost the same speed. If $f$ is compactly supported, then $T_\lambda\star f$ is even stationary. This is why to study the speed in the renewal theorem, we will study the speed in the convergence of equation~\ref{equation:renewal_theorem} and, actually, we will see in proposition~\ref{proposition:tempered_renewal_kernel} that this is more than a computational trick.

\medskip
\begin{lemma}\label{lemma:renewal_equation}
Let $\rho$ be a probability measure on $\R$ which have an exponential moment and a positive drift $\lambda$.

Then, for any non negative $f\in \cal S(\R)$ and any $x\in \R$,
\[
G\star f(x) = \lim_{s\to 0^+} \int_\R \widehat f(\xi) e^{i\xi x}\frac 1 {1-\widehat\rho(s-i\xi)}  \di \xi \text{ and }T_\lambda\star f(x) = \lim_{s\to 0^+} \frac{-1} \lambda \int_\R \widehat f(\xi) e^{i\xi x} \frac 1 {s-i\xi} \di \xi
\]
\end{lemma}

\begin{proof}
To prove this proposition, we are going to study the Fourier transform of $G \star \delta_x $. This is not well defined a priori so we need to approximate $\rho$ by another measure $\rho_s$. 

For $a\in\R$, note $\rho_a$ the measure on $\R$ defined by $\rho_a (A) = \int_A e^{-a y} \di\rho(y)$.
Since $\rho$ has exponential moments, the function $\varphi:a\mapsto \rho_a(\R)$ is defined on a neighbourhood of $0$ and is differentiable at $0$. Moreover, $\varphi(0)=1$ and $\varphi'(0)=-\int_\R y\di\rho(y) <0$ and so, there is $a_0\in\R_+^\star$ such that for any $a\in ]0,a_0[$, $\rho_a(\R)<1$.

Note $P_a$ the operator associated to $\rho_a$. Saying that $\rho_a(\R)<1$ means that for any bounded function $f$ on $\R$ and any $x\in \R$, $|P_a^nf(x)| \leqslant \|f\|_\infty \rho_a(\R)^n$ and so, $\sum_{n\in\N} P_a^n f$ is well defined.

\medskip
We are going to prove that for any non negative function $f\in\cal S(\R)$ and any $x\in \R$,
\[
\lim_{a\to 0^+} \sum_{n=0}^{+\infty} P_a^nf(x) = \sum_{n=0}^{+\infty} P^nf(x)
\]
and that for any $a\in]0,a_0[$,
\[
\sum_{n=0}^{+\infty} P_a^n f(x) = \int_\R \widehat f(\xi) e^{i\xi x} \frac 1 {1-\widehat\rho(a-i\xi)} \di \xi
\]

First,
\begin{flalign*}
\sum_{n=0}^{+\infty} P_a^nf(x) &= \sum_{n=0}^{+\infty} \int_\R e^{-ay} f(x+y) \di\rho^{\star n}(y) \\
&= \sum_{n=0}^{+\infty} \int_{\R_+^\star} e^{-ay} f(x+y) \di\rho^{\star n}(y) + \sum_{n=0}^{+\infty} \int_{\R_-} e^{-ay} f(x+y) \di\rho^{\star n}(y)
\end{flalign*}
And the monotone convergence theorem (we took $f$ non negative) proves that
\[
\lim_{a\to 0^+}  \sum_{n=0}^{+\infty} \int_{\R_+} e^{-ay} f(x+y) \di\rho^{\star n}(y) = \sum_{n=0}^{+\infty} \int_{\R_+}  f(x+y) \di\rho^{\star n}(y)
\]
Moreover,
\[
\int_{\R_-^\star} e^{-ay} f(x+y) \di\rho^{\star n}(y) \leqslant \|f\|_\infty \int_{\R_-^\star} e^{-ay} \di\rho^{\star n}(y) = \|f\|_\infty \rho_a(\R_+^\star)^n
\]
so, we also have that
\[
\lim_{a\to 0^+}  \sum_{n=0}^{+\infty} \int_{\R_-^\star} e^{ay} f(x+y) \di\rho^{\star n}(y) = \sum_{n=0}^{+\infty} \int_{\R_-^\star}  f(x+y) \di\rho^{\star n}(y)
\]
Thus,
\[
\lim_{a\to 0^+} \sum_{n=0}^{+\infty} P_a^nf(x) = \sum_{n=0}^{+\infty} P^nf(x)
\]

Finally, we took $f\in \cal S(R)$, this means that $\widehat f$, the Fourier-transform of $f$, also belong to $\cal S(\R)$.
So,
\begin{flalign*}
\sum_{n=0}^{+\infty} P_a^nf(x)& =  \sum_{n=0}^{+\infty} \int_\R f(x+y) e^{-ay}\di\rho^{\star n}(y) = \frac 1 {2\pi}\sum_{n=0}^{+\infty}\int_\R \int_\R \widehat f(\xi) e^{i\xi(x+y)} e^{-ay}\di\xi \di\rho^{\star n}(y)\\
&= \frac 1 {2\pi} \int_\R \sum_{n=0}^{+\infty}\int_\R \widehat f(\xi) e^{i\xi(x+y)} e^{-ay} \di\rho^{\star n}(y) \di\xi \\
&=\frac 1 {2\pi} \int_\R \widehat f(\xi) e^{i\xi x} \sum_{n=0}^{+\infty} \int_\R e^{-(a-i\xi)y} \di\rho^{\star n}(y) \di\xi \\
&= \frac 1 {2\pi} \int_\R \widehat f(\xi) e^{i\xi x} \sum_{n=0}^{+\infty} \widehat\rho(a-i\xi)^n \di\xi = \frac 1 {2\pi} \int_\R \widehat f(\xi) e^{i\xi x} \frac 1 {1- \widehat\rho(a-i\xi)} \di\xi
\end{flalign*}
The reader may remark that at any step, we can use Fubini's theorem as long as we don't change the order between the sum and the integral against $\rho$.

Finally,
\begin{flalign*}
\int_\R \widehat f(\xi) e^{i\xi x} \frac 1 {a-i\xi} \di \xi = \int_\R f(x+u) \un_{\R_+}(u) e^{-au} \di u \xrightarrow[a\to 0^+]\, \int_\R f(x+u)\un_{\R_+}(u) \di u
\end{flalign*}

So,
\[
\lim_{a\to 0^+} \frac{1} \lambda \int_\R \widehat f(\xi) e^{i\xi x} \frac 1 {a-i\xi} \di \xi = T_\lambda \star f(x)
\]
and this finishes the proof of the lemma.
\end{proof}

\subsection{Speed in the renewal theorem}

\begin{proposition}\label{proposition:tempered_renewal_kernel}
Let $\rho$ be a probability measure on $\R$ which have an exponential moment and a positive drift $\lambda$.

Then, $G-T_\lambda$ is a tempered distribution of fast decrease if and only if $\rho$ is weakly diophantine.

Moreover, in this case,
\[
\widehat{G-T_\lambda}(i\xi)= \frac 1 {1-\widehat \rho(-i\xi)} - \frac 1 \lambda \frac 1 {i\xi} 
\]
\end{proposition}

\begin{remark}
In particular, if $\rho$ is weakly-diophantine, then for any $f\in \cal S(\R)$, the convergence in the renewal theorem if faster than any polynomial.
\end{remark}

\begin{proof}
We showed in lemma~\ref{lemma:renewal_equation} that for any non negative $f\in \cal S(\R)$ and any $x\in \R$,
\[
(G-T_\lambda)\star f(x) = \lim_{s\to 0^+} \int_\R \widehat f(\xi) e^{i\xi x} \left(\frac 1 {1-\widehat \rho(s-i\xi)} + \frac 1 \lambda \frac 1 {s-i\xi} \right) \di \xi
\]
So, if $G-T_\lambda \in \cal O_c(\R)$, then
\[
\widehat{G-T_\lambda}(i\xi) = \frac 1 {1-\widehat \rho(-i\xi)} - \frac 1 \lambda \frac 1 {i\xi}
\]
and there is $l\in \N$ such that
\[
 \frac 1 {1-\widehat \rho(-i\xi)} - \frac 1 \lambda \frac 1 {i\xi} =\widehat{G-T_\lambda}(i\xi) \in \cal O(1+|\xi|^l)
\]
but, this exactly means that $\rho$ is $l-$weakly-diophantine.

On the other hand, if $\rho$ is $l-$weakly diophantine, then, we have that
\[
\frac 1 {1-\widehat \rho(-i\xi)} - \frac 1 \lambda \frac 1 {i\xi}  \in \cal O(1+|\xi|^l)
\]
and, as $\widehat f$ decreases faster than ant polynomial, the dominated convergence theorem proves that
\[
(G-T_\lambda)\star f(x) = \int_\R \widehat f(i\xi) e^{i\xi x} \left(\frac 1 {1-\widehat \rho(-i\xi)} - \frac 1 \lambda \frac 1 {i\xi} \right) \di \xi
\]
But, as $\rho$ is $l-$weakly diophantine, the function $z\mapsto \frac 1 {1-\widehat \rho(z)} + \frac 1 \lambda \frac 1 {z}$ is holomorphic on
\[
\left\{z\in \C \middle| |\Re(z)| \leqslant \frac C {1+|\Im(z)|^l} \right\}
\]
and satisfies that
\[
\sup_{t\in \R} \frac 1 {1+|t|^l} \left| \frac 1 {1-\widehat \rho(it)} + \frac 1 \lambda \frac 1 {it}\right| \text{ is finite}
\]
This means that it satisfies to the assumptions of the theorem 15 in \cite[Chapter 4, part 7]{GC64} and this proves that it is the Fourier-Laplace transform of a tempered distribution of
fast decrease.
\end{proof}

\begin{proposition}
Let $\rho$ be a probability measure on $\R$ which have an exponential moment and a negative drift $\lambda$.

Assume that $\rho$ is $l-$diophantine for some $l\in \R_+$.

Then, $G-T_\lambda$ is continuous from $\cal H^{ m+2+l,\infty}_{p_{k+2}}(\R)$ to $\cal H^{m,\infty}_{p_k}(\R)$.
\end{proposition}

\begin{proof}
The proof is direct from corollary~\ref{corollary:polynomial_conergence_norme_infinity}
\end{proof}

In the weakly diophantine case, it turns out that the speed cannot be faster than faster than any polynomial. We first show that it cannot be between polynomial and exponential and in proposition~\ref{proposition:obstruction_exponential} we will prove that it cannot be exponential either.
\begin{proposition}\label{proposition:not_intermediate_speed}
Let $\rho$ be a probability measure on $\R$ which have an exponential moment and a negative drift $\lambda$.

Assume that $\rho$ is weakly-diophantine but that there is $l\in \R_+^\star$ such that $\rho$ is \underline{not} $l-$weakly-diophantine.

Let $\omega\in \Omega$ (see definition~\ref{definition:Omega_0}) be such that $G-T_\lambda: \left\{\begin{array}{ccc} C^\infty_c(\R) & \to & \mathrm{L}^2_\omega(\R) \\ f & \mapsto & (G-T_\lambda)\star f\end{array}\right.$ is continuous, then $\omega$ doesn't grow faster than any polynomial.
\end{proposition}

\begin{proof}
We saw in lemma~\ref{lemma:renewal_equation} that for any $f\in \cal D(\R)$,
\[
(G-T_\lambda)\star f(x) = \int_\R \widehat f(\xi) e^{i\xi x} \left(\frac 1 {1-\widehat \rho(-i\xi)} + \frac 1 \lambda \frac 1 {-i\xi} \right) \di \xi
\]
For $z\in \C_\eta$ such that $\widehat \rho(z)\not=1$, note
\[
U (z) = \frac 1 {1-\widehat \rho(z)} + \frac 1 \lambda \frac 1 {z} 
\]
Then, $U$ is a meromorphic function on $\C_\eta$ whose poles are simple and uniformly isolated according to lemma~\ref{lemma:uniformily_isolated_zeros}.

Note $\cal A =\left\{z\in \C_\eta\setminus\{0\}\middle| \widehat\rho(z)=1\right\}$.
And
\[
U_0(z) = U(z) - \sum_{a\in \cal A} \left( \frac z a \right)^2 \frac 1 {\widehat \rho'(a)} \frac 1 {z-a}
\]
Then, $U_0$ is holomorphic on $\C_\eta$.

Moreover, using lemma~\ref{lemma:approximate_zeros}, we know that there are $\varepsilon,\delta\in \R_+^\star$ such that for any $z\in \C_\eta$, if
\[
\inf_{a\in \cal A \cup\{0\}} |a-z| \geqslant \varepsilon
\]
then,
\[
|1-\widehat \rho(z)| \geqslant \delta
\]
so, for any $z\in \C_\eta$ such that
\[
\inf_{a\in \cal A \cup\{0\}} |a-z| \geqslant \varepsilon
\]
we have that
\[
|U_0(z)| \leqslant \frac 1 \delta + \frac 1 \varepsilon \frac 1 {|\lambda|} + |z|^2 \frac 1 \varepsilon \sum_{a\in \cal A} \frac 1 {|a|^2} \frac 1 {|\widehat \rho'(a)|} 
\]

This means that
\[
\sup_{\substack{z\in \C_\eta \\ \di(z,\cal A\cup\{0\}) \geqslant \varepsilon}} \frac{|U_0(z)|}{1+|z|^2} \text{ is finite}
\]
But, for each $a\in \cal A$, we can apply the maximum principle on $B(a,2\varepsilon)$ (remind that $\cal A$ is a Delone set as we saw in remark~\ref{remark:Delone_set}) to $U_0/(B+z^2)$ for some large $B\in \R$, we finally get that
\[
\sup_{z\in \C_\eta}\frac{|U_0(z)|}{1+|z|^2} \text{ is finite}
\]
Moreover, as $\rho$ is not $l-$weakly-diophantine for some $l\in \R_+^\star$, lemma~\ref{lemma:zeroes_fourier_transform_only_if_part} proves that for any $\varepsilon\in \R_+^\star$, there is $(z_\varepsilon)\in \C_{\eta_0}$ such that
\[
\widehat\rho(z_\varepsilon) =1\text{ and }\Re(z_\varepsilon) \geqslant \frac{-\varepsilon}{1+|\Im(z_\varepsilon)|^l}
\]
and this means that $U\in \cal L_\eta$ (see lemma~\ref{lemma:other_definition_L_eta}) and we conclude with proposition~\ref{proposition:between_polynomial_exponential}.
\end{proof}

\subsection{Exponential speed of convergence}

\begin{proposition}\label{proposition:obstruction_exponential}
Let $\rho$ be a measure on $\R$ that have an exponential moment and a positive drift $\lambda = \int_\R y\di\rho(y)>0$.

Assume that there is $\gamma\in\R_+^\star$ such that for any $f\in \cal C^\infty_c(\R)$
\[
(G-T_\lambda)\star f(x) \in \cal O\left( e^{-\gamma |x|}\right)
\]
Then, $\rho$ is $0-$weakly-diophantine.
\end{proposition}

The idea of the proof is that if there is such a $\gamma$, then we can solve the equation $f=g-Pg$ for any $f\in \cal C^{\infty}_c(\R)$ such that $\int_\R f(x) \di x= 0$ with some function $g$ that vanishes exponentially fast at $\pm\infty$ but this is impossible if $\rho$ is not $0-$diophantine since in this case, there are many $P-$invariant measures in $\cal C^0_{\gamma}(\R)^\star$ as we saw in remark~\ref{remark:invariant_measures}.

\begin{proof}
Let $f\in \cal C^\infty_c(\R)$ be such that $\int_\R f(x)\di x=0$.

Note $g=\sum_{n=0}^{+\infty} P^nf= G\star f(x)$.

Then, $f=g-Pg$ and for $x\in\R$ large enough (actually, larger than $M$ if $\supp(f)\subset [-M,M]$), since $T_\lambda\star f$ is compactly supported, we have that
\[
g(x) \in \cal O\left(e^{-\gamma |x|} \right)
\]
Take now $z\in \C$ such that $\widehat \rho(z)=1$ and $|\Re(z)|<\gamma$ (this always exists since $\widehat \rho(0)=1$).

Then, $(x\mapsto g(x)e^{zx}) \in \mathrm{L}^1(\R)$ and so we can use Fubini's theorem in the following computations (we can always assume that $\gamma$ is small enough to also have that $(x\mapsto Pg(x)e^{zx}) \in \mathrm{L}^1(\R)$ since $\rho$ has an exponential moment)
\begin{flalign*}
\int_\R f(x) e^{zx} \di x &= \int_\R (g(x)-Pg(x)) e^{zx} \di x = \int_\R g(x) e^{zx} \di x - \int_{\R^2} g(x+y)  e^{zx} \di\rho(y) \di x \\
&= \int_\R g(x) e^{zx} \di x - \int_\R g(x)e^{zx} \di x \int_\R e^{-zy} \di\rho(y)\\
& = (1-\widehat\rho(z))\int_\R g(x) e^{zx} \di x = 0
\end{flalign*}
This means that for any $f\in \cal C^\infty_c(\R)$,
\[
\int_\R f(x)\di x =0 \Rightarrow \forall z\in \C_\gamma \text{ st } \widehat\rho(z)=1,\;\;\int_\R f(x) e^{zx} \di x=0
\]
Finally, take $f\in \cal C^\infty_c(\R)$ and note that $\int_\R f'(x)\di x=0$.

Then, for any $z\in \C_\gamma$ such that $\widehat\rho(z)=1$,
\[
0= \int_\R f'(x) e^{zx} \di x = z \int_\R f(x) e^{zx} \di x
\]
Thus, for any $f\in \cal C^\infty_c(\R)$ and any $z\in \C_\gamma$ such that $\widehat\rho(z)=1$,
\[
z\int_\R f(x) e^{zx}\di x = 0
\]
And so,
\[
\left\{z\in \C_\gamma \middle|\widehat\rho(z)=1 \right\} = \{0\}
\]
And finally, proposition~\ref{proposition:zeroes_fourier_transform} proves that $\rho$ is $0-$weakly-diophantine.
\end{proof}

\begin{proposition}\label{proposition:renewal_theorem_exponential}
Let $\rho$ be a measure on $\R$ that have an exponential moment and a negative drift $\lambda = \int_\R y\di\rho(y)<0
$.

If $\rho$ is $0-$weakly-diophantine then there is $\gamma_0\in\R_+^\star$ such that for any $\gamma\in ]0,\gamma_0[$, the operator $G-T_\lambda$ is a continuous Fourier multiplier from $\cal C^{2}_\gamma(\R)$ to $\cal C^{0}_\gamma(\R)$ and from $\{f\in \mathrm{L}^2(\R)| e^{\gamma|x|} f(x) \in \mathrm{L}^2(\R)\}$ to itself.
\end{proposition}

\begin{remark}
In particular, if $\rho$ is $0-$weakly-diophantine, then for any $f\in \cal C^2_\gamma$, the speed in the renewal theorem is exponential.
\end{remark}

\begin{proof}
Let $f\in \cal S(\R)$, then,
\begin{flalign*}
(G-T_\lambda) \star f(x) &= \sum_{n=0}^{+\infty } P^n f(x) +\frac 1 \lambda \int_{x} ^{+\infty} f(t)\di t \\
& =\lim_{s\to 0^+} \int_\R \widehat f(i\xi) e^{i\xi x}\left( \frac 1 {1-\widehat \rho(s-i\xi)} -\frac 1 \lambda \frac 1 {s+i\xi} \right)\di \xi \\&= \int_\R \widehat f(i\xi) e^{i\xi x} \left( \frac 1 {1-\widehat \rho(-i\xi)} - \frac 1 \lambda \frac 1 {i\xi} \right) \di \xi
\end{flalign*}
Moreover $z\mapsto \frac 1 {1-\widehat \rho(-z)} - \frac 1 \lambda \frac 1 z$ is holomorphic and bounded on $\C_{\gamma_0}$ for some $\gamma_0\in\R_+^\star$.

So, we can apply lemma~\ref{lemma:reed_simon_exponential_speed} with $M=2$ and corollary~\ref{corollary:exponential_sobolev} with $l=0$.
\end{proof}

\begin{remark}
In proposition~\ref{proposition:renewal_theorem_exponential}, we have to take functions in $\cal C^2(\R)$. The issue is the same as the one discussed in remark~\ref{remark:too_restrictive_assumption}
\end{remark}

\appendix
\section{Dirac combs on strips in the complex plane}

Next proposition is a generalization of the fact that a weighted Dirac comb on $\R$ is a tempered distribution if and only if it's weight is polynomial.

\medskip
For $z\in \C$, we note $\widehat \delta_z$ the distribution that maps $\varphi\in \cal C^\infty_c(\R)$ onto $\widehat \varphi(z)$.

\begin{proposition}\label{proposition:strip_control_Laplace_dirac}
Let $\eta\in \R_+^\star$, $\cal A \subset \C_\eta$ and $(\lambda_a)\in \C^{\cal A}$ be such that for any $A\in \R_+^\star$,
\[
\sup_{\substack{a\in \cal A\\|a| \leqslant A }} |\lambda_a| \text{ is finite}
\]

Then, $(\lambda_a \widehat \delta_a)_{a\in \cal A}$ is bounded in $\cal D'(\R)$ if and only if there is $m\in \N$ such that $(\lambda_a) \in \cal O(1+|a|^m)$.
\end{proposition}

\begin{proof}
First, assume that there are $C\in \R_+$ and $m\in \N$ such that for any $a\in \cal A$,
\[
|\lambda_a| \leqslant C(1+|a|^m)
\]
Let $M\in\R_+^\star$, $f\in \cal C^\infty([-M,M])$ and $z\in \C_\eta$.
\[
\widehat f(z)=\int_\R e^{-zx} f(x) \di x = \frac 1 {z^n} \int_\R e^{-zx} f^{(n)}(x) \di x
\]
thus,
\[
|z|^n |\widehat f(z)| \leqslant \int_\R e^{-\Re(z)x} |f^{(n)}(x)| \di x \leqslant  2Me^{\eta M} \|f^{(n)}\|_\infty
\]
so, for any $a\in \cal A$,
\[
|\lambda_a| |\widehat f(a)| \leqslant C(1+|a|)^m  |\widehat f(a)|\leqslant 4CMe^{\eta M} \max_{k\in \{0,\dots, m\}} \|f^{(k)}\|_\infty
\]
and this proves that $(\lambda_a\widehat \delta_a)_{a\in \cal A}$ is bounded in $\cal D'(\R)$.

\medskip
To prove the converse, assume now that $(\lambda_a\widehat \delta_a)_{a\in \cal A}$ is bounded in $\cal D'(\R)$.

Let $f\in \cal C^\infty_c(\R)$ be such that $\widehat f(i)\not=0$.

For $t\in \R_+^\star$, note $f_t(x) = tf(tx)$.

Then, for any $z\in \C$,
\[
\widehat f_t(z) =\int_\R e^{-zx} t f(tx) \di x=\int_\R e^{-zu/t} f(u)\di u= \widehat f(\frac z t)
\]

Take $M\in \R_+^\star$ such that $\supp f\subset [-M,M]$. Then, we also have that for any $t\in [1,+\infty[$, $\supp(f_t) \subset [-M,M]$.

Saying that $(\lambda_a \widehat \delta_a)$ is bounded in $\cal D'(\R)$ means that there are $C,N$, depending on $M$ such that for any $\varphi\in \cal C^\infty([-M,M])$,
\[
\sup_{a\in \cal A} |\lambda_a| |\widehat \varphi(a)| \leqslant C \max_{k\in [ 0,N ]} \|\varphi^{(k)}\|_\infty 
\]
In particular, we get that for any $a\in \cal A$ such that $|a| \geqslant 1$,
\[
|\lambda_a| |\widehat f_{|a|} (a)| \leqslant C \max_{k\in [ 0,N ]} \|f_{|a|}^{(k)}\|_\infty 
\]
Remark that for any $t\in \R_+^\star$ and any $k\in \N$, $\|f_t^{(k)}\|_\infty \leqslant t^{k+1} \|f^{(k)}\|_\infty$.

So, for any $a\in \cal A$ such that $|a| \geqslant 1$,
\[
|\lambda_a| \left| \widehat f\left( \frac a {|a|} \right) \right| \leqslant C \left(1+|a|^{N+1}\right) \max_{k\in [ 0,N ]} \|f^{(k)}\|_\infty
\]
Therefore, the proof will be finished if we manage to prove that there is $A\in [1,+\infty[$ such that
\begin{equation}\label{equation:control_comb}
\inf_{a\in \cal A, |a| \geqslant A} \left|\widehat f\left(\frac a{|a|}\right)\right| >0
\end{equation}
since by assumption, we also have that
\[
\sup_{a\in \cal A, |a|\leqslant A} |\lambda_a| \text{ is finite}
\]
For $z\in \C_\eta$, such that $\Im(z) \geqslant 0$, we have that
\[
\left| z - i|z| \right| \leqslant |\Re(z)| + \left|\Im(z)-|z|\right| \leqslant 2|\Re(z)| \leqslant 2\eta
\]
So, for any $z\in \C_\eta\setminus\{0\}$ such that $\Re(z)\geqslant 0$,
\[
\left|\frac{z}{|z|}-i\right| \leqslant \frac{2\eta}{|z|}
\]
in the same way, we get that for any $z\in \C_\eta$ such that $\Im(z)\leqslant 0$,
\[
\left|\frac{z}{|z|}+i\right| \leqslant \frac{2\eta}{|z|}
\]

Moreover, for any $z\in \C_\eta$, $\left|\widehat f(\overline z)\right|$ so $\left|\widehat f(i) \right|=\left|\widehat f(-i)\right|\not=0$.

What we just proved (as $\widehat f$ is continuous) is that for any $\varepsilon \in \R_+^\star$, there is $A_\varepsilon \in \R_+$ such that for any $z\in \C_\eta$ such that $|z| \geqslant A_\varepsilon$,
\[
\left| \widehat f\left( \frac z{|z|} \right) - \widehat f\left(sign(\Im(z))i \right) \right| \leqslant \varepsilon
\]
and as $\widehat f(i) \not=0$, and as we may assume without any loss of generality that $A_\varepsilon\geqslant 1$, this is exactly what we wanted to prove in equation~\ref{equation:control_comb}.
\end{proof}

\section{Fourier-Laplace transforms of non integrable functions} \label{appendix:Omega}

In this section, we note $\C_+ = \{z\in \Z| \Re(z) \geqslant 0 \}$ and $\C_+^\star = \{z\in \C; \Re(z) >0\}$.

Let $\omega$ be a function on $\R_+$ such that for any $\varepsilon \in \R_+^\star$,
\begin{equation}\label{equation:omega_sous_expo}
\lim_{x\to +\infty} e^{-\varepsilon x} \omega(x) = 0
\end{equation}
We call such functions strictly sub-exponential.

\medskip
For $z\in\C_+^\star$, the one-sided Laplace transform of $\omega$ is well defined at $z$ :
\[
\cal L(\omega)(z) = \int_{\R_+} e^{-zx} \omega(x) \di x
\]
Moreover, if $\omega \in \mathrm{L}^1(\R)$, then Jensen's inequality proves that the integral is actually bounded on $\C_+^\star$.

If $\omega \not\in \mathrm L^1(\R)$ and is non negative, the monotone convergence theorem proves that
\[
\lim_{z\to 0^+,\;z\in \R} \int_{\R_+} \omega(x) e^{-zx}\di x= \int_{\R_+} \omega(x) \di x = +\infty
\]

On the contrary, if we fix $a\in \R_+^\star$, $\left(x\mapsto e^{-ax}\omega(x) \right) \in \mathrm{L}^1(\R_+)$ (because of \ref{equation:omega_sous_expo}) so Riemann-Lebesgue's theorem proves that
\[
\lim_{b\to \pm\infty,\;b\in \R} \int_{\R_+} \omega(x)e^{-(a+ib)x} \di x = 0
\]

In this section, we want to understand a little bit more this phenomenon.

More specifically, we want to find functions $\omega$ on $\R_+$ such that, if we note
\[
\cal U_{\delta} = \left\{z\in \C_+^\star\middle| |z|\geqslant \delta \right\}
\]
we have that for any $\delta\in \R_+^\star$,
\begin{equation} \label{equation:theta_omega}
\Theta_\omega(\delta) = \sup_{z\in \cal U_{\delta}} \left| \int_{\R_+} \omega(x) e^{-zx} \di x \right|
\end{equation}
is finite. Or in other words, we want to find strictly sub-exponentials functions whose one-sided Fourier-Laplace transform is bounded on $\C_+^\star$ except at $0$.

Thus, we make the following
\begin{definition}\label{definition:Omega_0}
For a borelian function $\omega$ on $\R_+$ let $\Theta_\omega$ be the function on $\R_+$ defined by equation~\ref{equation:theta_omega}.

Then, we note
\[
\Omega_0 = \left\{\omega \in \cal C^0(\R_+)\middle| \omega \text{ is strictly subexponential and } \forall \delta\in\R_+^\star\; \Theta_\omega(\delta) \text{ is finite} \right\}
\]
and
\[
\Omega = \left\{\begin{array}{c|l} \multirow{3}{*}{$\omega \in \cal C^0(\R, [1,+\infty[) $}& \omega \text{ is even} \\ &\forall \varepsilon\in \R_+^\star\; \lim_{x\to \pm\infty} e^{-\varepsilon|x|} \omega(x) = 0 \\ &\forall \delta \in \R_+^\star\; \sup_{\substack{z\in \C\\ \Re(z)>0 \text{ and }|z| \geqslant \delta}} \left| \int_{0}^{+\infty} \omega(x) e^{-zx} \di x \right| \text{ is finite}\end{array} \right\}
\]
\end{definition}

It is not clear at all to find non trivial functions belonging to $\Omega_0$ (i.e. functions that are not in $\mathrm{L}^1(\R)$) and the aim of this section is to prove proposition~\ref{proposition:Laplace_transform_analytic_functions} which gives many examples and lemma~\ref{lemma:essential_singularity_Laplace_transform} which characterizes those of our examples that are bounded by polynomials.

\medskip
But first, we gathered the easy examples in next
\begin{lemma}~
\begin{enumerate}[label=(\roman*)]
\item\label{item:omega_vector} The set $\Omega_0$ is a vector space and moreover, if $\omega\in \Omega_0$, any function that is equal to $\omega$ out of a bounded set is in $\Omega_0$.
\item\label{item:omega_integration} The set $\Omega_0$ is stable by integration : it means that if $\omega\in \Omega_0$ then any primitive of $\omega$ also belong to $\Omega_0$.
\item\label{item:omega_constants} The constant functions belong to $\Omega_0$ (and so $\Omega_0$ also contains the polynomials)
\item\label{item:omega_delay} If $\omega\in \Omega_0$, then for any $t\in \R_+^\star$, $(x\mapsto \omega(x+t))\in \Omega_0$.
\end{enumerate}
\end{lemma}

\begin{proof}
The points \ref{item:omega_vector}, \ref{item:omega_constants} and \ref{item:omega_delay} are direct.

To prove \ref{item:omega_integration}, let $f\in \Omega_0$. And note $F$ a primitive of $\omega$.

Then, for any $\varepsilon \in \R_+^\star$, $M=\frac 1 \varepsilon\sup_{x\in\R_+} e^{-\varepsilon x} |f(x)|$ is finite by definition of $\Omega_0$ and for any $x\in \R_+$,
\begin{flalign*}
e^{-\varepsilon x} |F(x)-F(x-1)| &= e^{-\varepsilon x} \left| \int_{x-1}^x f(t) \di t \right|\leqslant \sup_{y\in\R_+} e^{-\varepsilon y} |f(y)| \int_{x-1}^x e^{-\varepsilon x} e^{\varepsilon t} \di t \leqslant M
\end{flalign*}
and
\begin{flalign*}
e^{-\varepsilon x} |F(x)| &\leqslant e^{-\varepsilon x} \left| F(x) - F(x-1)\right| + e^{-\varepsilon} e^{-\varepsilon(x-1)} |F(x-1)| \\
& \leqslant M + e^{-\varepsilon} \left(M+ e^{-\varepsilon}e^{-\varepsilon(x-2)} |F(x-2)|\right) \\
&\leqslant \frac M {1-e^{-\varepsilon}} + e^{-\varepsilon\lfloor x \rfloor} \sup_{t\in [0,1]} |F(t)|
\end{flalign*}

So, for any $\varepsilon \in \R_+^\star$,
\[
\lim_{x\to +\infty} e^{-\varepsilon x} |F(x)| =0
\]
And this means that the Laplace transform of $F$ is well definite on $\C_+^\star$.

Moreover for any $z\in \C_+^\star$,
\[
\cal L (F) (z) = \int_{\R_+} e^{-zx} F(x) \di x = \left[ \frac{e^{-zx}}{-z} F(x)\right]_0^{+\infty} - \int_{\R_+} \frac{e^{-zx}}{-z} f(x) \di x = \frac 1 z F(0) + \frac 1 z \cal L(f)(z)
\]
Hence, $F\in \Omega_0$.
\end{proof}

A first non easy example of a function belonging to $\Omega$ is given by next \begin{lemma}\label{lemma:Laplace_transform_sqrt_function}
For any $A\in \R$ and any $\alpha\in]0,1[$, $(x\mapsto e^{Ax^\alpha})\in \Omega$.
\end{lemma}

\begin{proof}
We are going to prove, as a first step, that for any $A\in \R_+^\star$ and $\alpha\in ]0,1[$ and any $z\in \C_+^\star$
\begin{equation}\label{equation:control_subexponential_functions}
\int_{\R_+} e^{Ax^\alpha - zx} \di x =\frac {1} z \sum_{k=0}^{+\infty}\left( \frac A  {z^\alpha} \right)^k \frac {\Gamma(\alpha k+1)} {\Gamma(k+1)}   
\end{equation}
To do so, we are first going to prove that the formula holds for $z\in \R_+^\star$.

Indeed, we have, for $t\in \R_+^\star$,
\begin{flalign*}
\int_{\R_+} e^{Ax^\alpha -tx } \di x& = \sum_{k\in \N} \frac{A^k}{k!} \int_{\R_+} x^{\alpha k} e^{-tx} \di x = \sum_{k\in\N}  \frac{A^k}{k!} \frac 1 {t^{\alpha k+1}} \int_{\R_+} u^{\alpha k} e^{-u} \di u \\&= \frac 1 t  \sum_{k=0}^{+\infty}\left( \frac A  {t^\alpha} \right)^k \frac {\Gamma(\alpha k+1)} {\Gamma(k+1)}   
\end{flalign*}
and this proves formula~\ref{equation:control_subexponential_functions} for any $z\in \C_+^\star$ since $\Gamma(k+1) = k!$.

\medskip
As the left side of the equation is an analytic function of $z$, we only need to prove that so does the right side (we have fixed here a determination of the logarithm on $\C_+^\star$ to define $z^\alpha$) and the isolated zero theorem will then show that the formula holds for any $z\in \C_+^\star$.

Using Jensen's inequality, we get that for any $z\in \C_+^\star$,
\[
\left|\frac {1} z \sum_{k=0}^{+\infty}\left( \frac A  {z^\alpha} \right)^k \frac {\Gamma(\alpha k+1)} {\Gamma(k+1)} \right| \leqslant \frac 1 {|z|} \sum_{k=0}^{+\infty} \left| \frac A  {|z|^\alpha} \right|^k \frac {\Gamma(\alpha k+1)} {\Gamma(k+1)}   
\]
Therefore, for any $z\in \cal U_\delta$,
\[
\left|\frac {1} z \sum_{k=0}^{+\infty}\left( \frac A  {z^\alpha} \right)^k \frac {\Gamma(\alpha k+1)} {\Gamma(k+1)} \right| \leqslant \frac 1 {\delta} \sum_{k=0}^{+\infty} \left| \frac A  {\delta^\alpha} \right|^k \frac {\Gamma(\alpha k+1)} {\Gamma(k+1)}
\]
And this last sum is finite for any $\delta\in \R_+^\star$ since Stirling's formula proves that
\[
\frac{\Gamma(\alpha k+1)}{\Gamma(k+1)} \sim \frac{\sqrt{k+1}}{\sqrt{\alpha k +1}} e^{k(1-\alpha)} \frac{ (\alpha k+1)^{\alpha k +1}}{(k+1)^{k+1}} \in \cal O(k^{-k(1-\alpha')}) \text{ for any }\alpha'\in ]\alpha,1[
\]
this proves both that the right side of equation~\ref{equation:control_subexponential_functions} is analytic and our lemma.
\end{proof}

Remark that the kernel $(x\mapsto e^{-zx})$ is not non negative so it is not clear at all that if $\omega\in \Omega_0$ and if $\tilde \omega$ is such that $0\leqslant \tilde \omega \leqslant \omega$ then $\tilde \omega\in \Omega_0$. The problem is that we can't use Jensen's inequality to bound the Laplace transform of $\tilde \omega$ at some point $z$ since it would kill the imaginary part of $z$. However, we can adapt in some way the proof of the previous lemma to ``remove the z from the integral'' before applying Jensen's inequality.

This is what we do in next proposition, but the price to pay with this method is that we are only able to deal with functions that have an holomorphic continuation on $\C_+^\star$.

\begin{proposition}\label{proposition:Laplace_transform_analytic_functions}
Let $\Phi$ be an analytic function on $\C_+^\star$ such that for some $\alpha\in ]0,1[$,
\[
\sup_{z\in \C_+^\star} \frac{|\Phi(z)|}{|1+z|^\alpha} \text{ is finite}
\]
Then, for any $\delta\in \R_+^\star$,
\[
\sup_{z\in \cal U_\delta} \left| \int_{\R_+} e^{\Phi(x)} e^{-zx} \di x \right| \text{ is finite}
\]
in particular, if $\Phi(\R_+)\subset \R_+$, then the function $\left(x\mapsto e^{\Phi(|x|)}\right)$ belongs to $\Omega$.
\end{proposition}

\begin{remark}
Example of such functions are $\Phi(x) = A(\ln(1+x))^M$ for $A,M\in \R$, or any sum or composition of such functions as long as at any step, they map $\{z\in \C|\Re(z)>0\}$ on to itself (e.g. if $l(x) = \ln(1+x)$, then for any $m\in \N^\star$, $l$ composed $m$ times works).
\end{remark}

\begin{remark}
For $M=1$, we have that $e^{A(\ln(1+|x|)) ^M} = (1+|x|)^A$ and we find that the polynomial functions belong to $\Omega$.
\end{remark}

\begin{proof}
The proof is almost the same as the one of lemma~\ref{lemma:Laplace_transform_sqrt_function}.

For $z\in \C_+^\star$, note
\[
\Psi(z) = \frac{\Phi(z)}{(1+z)^\alpha}
\]

Then, $\Psi$ is analytic and bounded on $\C_+^\star$ by assumption.

Moreover, we claim that for any $z\in \C_+^\star$, we have that
\[
\int_{\R_+} e^{\Phi(x)} e^{-zx} \di x = \frac {1} z \sum_{k=0}^{+\infty} \frac 1 {k!} \int_{\R_+} \left(1+ \frac u z\right)^{k\alpha} \left(\Psi\left(\frac u {z} \right)\right)^k e^{-u} \di u
\]
We use the same technique since both sides are analytic function of $z$ and they are equal on $\R_+^\star$.

And this proves that for any $z\in \C_+^\star$,
\[
\left|\int_{\R_+} e^{\Phi(x)} e^{-zx} \di x\right| \leqslant \frac 1 {|z|} \sum_{k=0}^{+\infty} \frac {A^k} {k!}  \int_{\R_+} \left| 1 + \frac u z \right|^{k\alpha} e^{-u} \di u   \text{ where }A=\|\Psi\|_\infty
\]
But, $\Re(z)>0$ and $u\in \R_+$, so
\[
\left| 1 + \frac u z \right|^{k\alpha} \leqslant \left| 1 + \frac u z \right|^{\lceil k\alpha\rceil} = \sum_{l=0}^{\lceil k\alpha\rceil} \binom {\lceil k\alpha\rceil} l \left|\frac u z\right| ^{\lceil k\alpha\rceil-l}
\]
Therefore,
\[
\left|\int_{\R_+} e^{\Phi(x)} e^{-zx} \di x\right| \leqslant \frac 1 {|z|} \sum_{k=0}^{+\infty} \frac {A^k} {k!}  \sum_{l=0}^{\lceil k\alpha\rceil}  \binom {\lceil k\alpha\rceil} l |z|^{l-\lceil k\alpha\rceil} \int_{\R_+}u^{\lceil k\alpha\rceil-l}  e^{-u} \di u
\]
so,
\begin{flalign*}
\left|\int_{\R_+} e^{\Phi(x)} e^{-zx} \di x\right| &\leqslant \frac 1 {|z|} \sum_{k=0}^{+\infty} \frac {A^k} {k!}  \sum_{l=0}^{\lceil k\alpha\rceil}  \binom {\lceil k\alpha\rceil} l |z|^{l-\lceil k\alpha\rceil} \Gamma(\lceil k\alpha\rceil-l+1) \\
&= \frac 1 {|z|} \sum_{k=0}^{+\infty} {A^k} \sum_{l=0}^{\lceil k\alpha\rceil}  \frac{\lceil k\alpha \rceil !}{k! l!}  |z|^{l-\lceil k\alpha\rceil}\\
& \leqslant \frac 1 {|z|} \sum_{k=0}^{+\infty} \frac {A^k} {|z|^{k\alpha}} \frac{(\lceil k\alpha\rceil)!}{k!} \max( 1, |z|^{k\alpha}) e^1 \\
&\retrait\retrait \text{ since }|z|^l \leqslant  \max(1,|z|^l) \leqslant \max(1,|z|^{\lceil k\alpha\rceil}) \\
& \leqslant \frac e {|z|} \sum_{k=0}^{+\infty} \frac {A^k}{\min(1,|z|) ^{\lceil k\alpha\rceil}} \frac{(\lceil k\alpha\rceil)!}{k!} \leqslant  \frac e {|z|} \sum_{k=0}^{+\infty} \left(\frac {A}{\min(1,|z|)} \right)^k \frac{(\lceil k\alpha\rceil)!}{k!}
\end{flalign*}

Therefore, assuming without any loss of generality that $\delta \leqslant 1$, we have that for any $z\in \cal U_\delta$,
\[
\left|\int_{\R_+} e^{\Phi(x) -zx } \di x \right| \leqslant \frac e {\delta} \sum_{k=0}^{+\infty} \left( \frac A  {\delta} \right)^k \frac {\Gamma(\lceil\alpha k\rceil+1)} {\Gamma(k+1)} < +\infty
\]
where the finiteness of the sum comes from the fact that, as proved by Stirling's formula, for any $\alpha'\in]\alpha,1[$,
\[
 \frac{(\lceil k\alpha\rceil)!}{k!} \in \cal O ( k^{-k(1-\alpha')} )
\]
\end{proof}

In next lemma, we characterize the elements of $\Omega$ that are bounded by some polynomial. This will be useful to prove that the speed in the renewal theorem is not faster than faster than any polynomial.

\begin{lemma}\label{lemma:essential_singularity_Laplace_transform}Let $\omega : \R_+\to \R_+$ be a continuous function that doesn't vanish and such that for any $\varepsilon \in \R_+^\star$, $\lim_{x\to +\infty} e^{-\varepsilon x} \omega(x) =0$.

Assume that for some $l\in \N^\star$,
\[
x^l \in \cal O(\omega(x))
\]
Then,
\[
\liminf_{s\to 0^+} s^l \int_{\R_+} \omega(x) e^{-sx} \di x =+\infty
\]
\end{lemma}

\begin{proof}
Let $l\in \N$ and assume that $x^l \in \cal O(\omega(x))$. This means that there are $C,x_0$ such that for any $x\in [x_0,+\infty[$, $x^l \leqslant C_0 \omega(x)$. And, as $\omega$ is continuous and doesn't vanish on $\R_+$, the function $(x\mapsto x^l/\omega(x))$ is continuous on $\R_+$ so it is bounded on $[0,x_0]$ thus, there is $C\in \R_+^\star$ such that for any $x\in \R_+$, 
\[
x^l \leqslant C \omega(x)
\]
So, for any $s\in \R_+^\star$,
\[
\int_{\R_+} \omega(x) e^{-sx} \di x \geqslant \frac 1 {C}\int_{\R_+} x^{l} e^{-sx} \di x = \frac 1 {C} \frac 1 {s^{l+1}} \Gamma(l+1)
\]
and so, multiplying each side of the previous inequality by $s^l$ and letting $s$ going to $0^+$, we get the expected result.
\end{proof}

\bibliographystyle{amsalpha}
\bibliography{biblio}

\providecommand{\bysame}{\leavevmode\hbox to3em{\hrulefill}\thinspace}
\providecommand{\MR}{\relax\ifhmode\unskip\space\fi MR }
% \MRhref is called by the amsart/book/proc definition of \MR.
\providecommand{\MRhref}[2]{%
  \href{http://www.ams.org/mathscinet-getitem?mr=#1}{#2}
}
\providecommand{\href}[2]{#2}
\begin{thebibliography}{Gam01}

\bibitem[BG07]{BG07}
J.~Blanchet and P.~Glynn, \emph{Uniform renewal theory with applications to
  expansions of random geometric sums}, Adv. in Appl. Probab. \textbf{39}
  (2007), no.~4, 1070--1097. \MR{2381589 (2009e:60191)}

\bibitem[Bou04]{Bou04}
Tahar~Zam{\`e}ne Boulmezaoud, \emph{On the invariance of weighted {S}obolev
  spaces under {F}ourier transform}, C. R. Math. Acad. Sci. Paris \textbf{339}
  (2004), no.~12, 861--866. \MR{2111723 (2005g:42052)}

\bibitem[Bre05]{Bre05}
E.~Breuillard, \emph{Distributions diophantiennes et th\'eor\`eme limite local
  sur {$\Bbb R^d$}}, Probab. Theory Related Fields \textbf{132} (2005), no.~1,
  39--73. \MR{2136866 (2006b:60093)}

\bibitem[Car83]{Car83}
Hasse Carlsson, \emph{Remainder term estimates of the renewal function}, Ann.
  Probab. \textbf{11} (1983), no.~1, 143--157. \MR{682805 (84e:60127)}

\bibitem[Fel71]{Fel71}
William Feller, \emph{An introduction to probability theory and its
  applications. {V}ol. {II}.}, Second edition, John Wiley \& Sons, Inc., New
  York-London-Sydney, 1971. \MR{0270403 (42 \#5292)}

\bibitem[Gam01]{Gam01}
T.W. Gamelin, \emph{Complex analysis}, Undergraduate Texts in Mathematics,
  Springer, 2001.

\bibitem[GC64]{GC64}
I.~M. Guelfand and G.~E. Chilov, \emph{Les distributions. {T}ome 2: {E}spaces
  fondamentaux}, Traduit par Serge Vasilach. Collection Universitaire de
  Math\'ematiques, XV, Dunod, Paris, 1964. \MR{0161145 (28 \#4354)}

\bibitem[Har10]{Har09}
Zen Harper, \emph{Laplace transform representations and {P}aley-{W}iener
  theorems for functions on vertical strips}, Doc. Math. \textbf{15} (2010),
  235--254. \MR{2628844 (2011m:30077)}

\bibitem[L{\"o}f83]{Lof83}
J{\"o}rgen L{\"o}fstr{\"o}m, \emph{A nonexistence theorem for translation
  invariant operators on weighted {$L_{p}$}-spaces}, Math. Scand. \textbf{53}
  (1983), no.~1, 88--96. \MR{733941 (85d:42018)}

\bibitem[RS75]{RS75}
Michael Reed and Barry Simon, \emph{Methods of modern mathematical physics.
  {II}. {F}ourier analysis, self-adjointness}, Academic Press [Harcourt Brace
  Jovanovich, Publishers], New York-London, 1975. \MR{0493420 (58 \#12429b)}

\bibitem[Rud91]{Rud91}
Walter Rudin, \emph{Functional analysis}, second ed., International Series in
  Pure and Applied Mathematics, McGraw-Hill, Inc., New York, 1991. \MR{1157815
  (92k:46001)}

\bibitem[Smi54]{Smi64}
Walter~L. Smith, \emph{Asymptotic renewal theorems}, Proc. Roy. Soc. Edinburgh.
  Sect. A. \textbf{64} (1954), 9--48. \MR{0060755 (15,722f)}

\end{thebibliography}
\end{document}